\documentclass{amsart}

\newtheorem{theorem}{Theorem}[section]
\newtheorem{lemma}{Lemma}[section]
\newtheorem{proposition}{Proposition}[section]

\theoremstyle{definition}

\newtheorem{example}{Example}[section]

\newtheorem{corollary}{Corollary}[section]

\theoremstyle{remark}
\newtheorem{remark}{Remark}[section]

\numberwithin{equation}{section}

\usepackage[shortlabels]{enumitem}
\setlist[enumerate,1]{label={\upshape(\roman*)}}

\newcommand*\diff{\mathop{}\!\mathrm{d}}

\def\Res{\mathop{\text{Res}}\displaylimits}

\renewcommand\leftmark{L. LAMGOUNI}
\renewcommand\rightmark{C-POLYNOMIALS AND LC-FUNCTIONS}

\usepackage{textcmds}

\usepackage{hyperref}

\usepackage{cleveref}

\usepackage{textcmds}

\newcommand{\floor}[1]{\left\lfloor #1 \right\rfloor}

\usepackage{graphicx}
\graphicspath{ {./Figs/} }
\usepackage{caption}
\usepackage{float}

\begin{document}

\title{C-polynomials and LC-functions: towards a generalization of the Hurwitz zeta function}

\author{Lahcen Lamgouni}
\curraddr{Minister of National Education Preschool and Sports: Errachidia, Morocco}
\email{lahcen.lamgouni@gmail.com}






\begin{abstract}
Let $f(t)=\sum_{n=0}^{+\infty}\frac{C_{f,n}}{n!}t^n$ be an analytic function at $0$, and let $C_{f, n}(x)=\sum_{k=0}^{n}\binom{n}{k}C_{f,k} x^{n-k}$ be the sequence of Appell polynomials, referred to as \textit{C-polynomials associated to} $f$, constructed from the sequence of coefficients $C_{f,n}$. We also define $P_{f,n}(x)$ as the sequence of C-polynomials associated to the function $p_{f}(t)=f(t)(e^t-1)/t$, called \textit{P-polynomials associated to} $f$. This work investigates three main topics. Firstly, we examine the properties of C-polynomials and P-polynomials and the underlying features that connect them. Secondly, drawing inspiration from the definition of P-polynomials and subject to an additional condition on $f$, we introduce and study the complex-variable function $P_{f}(s,z)=\sum_{k=0}^{+\infty}\binom{z}{k}P_{f,k}s^{z-k}$, which generalizes the $s^z$ function and is denoted by $s^{(z,f)}$. Thirdly, the paper's significant contribution is the generalization of the Hurwitz zeta function and its fundamental properties, most notably Hurwitz's formula, by constructing a novel class of functions defined by $L(z,f)=\sum_{n=n_{f}}^{+\infty}n^{(-z,f)}$, which are intrinsically linked to C-polynomials and referred to as \textit{LC-functions associated to} $f$ (the constant $n_{f}$ is a positive integer dependent on the choice of $f$). This research offers a detailed analysis of C-polynomials, P-polynomials, and LC-functions associated to a given analytic function $f$, thoroughly examining their interrelations and introducing unexplored research directions for a novel and expansive class of LC-functions possessing a functional equation equivalent to that of the Riemann zeta function, thereby highlighting the potential applications and implications of the findings.
\end{abstract}

\maketitle

\section{Introduction}\label{sec1}

The Bernoulli polynomials hold considerable importance in numerous areas of pure and applied mathematics due to their diverse applications. These polynomials are derived from the sequence of Bernoulli numbers, represented by $B_{n}$, and can be expressed as $B_{n}(x)=\sum_{k=0}^{n}\binom{n}{k}B_{k} x^{n-k}$ (for a comprehensive understanding of their properties, refer to \cite{Arakawa} and Chapter 24 of \cite{Olver}). The prominence of Bernoulli polynomials can be attributed to their presence in the values of the Hurwitz zeta function $\zeta(z,a)$ at non-positive integers\footnote{The Hurwitz zeta function $\zeta(z,a)$ is formulated for $\Re(z)>1$ through the series $\zeta(z,a)=\sum_{n=0}^{+\infty}(n+a)^{-z}$, wherein $0<a\leq1$ signifies a constant real number. For every integer $n\geq0$, $\zeta(-n,a)=-\frac{B_{n+1}(a)}{n+1}$.}. It should be noted that the Hurwitz zeta function make it possible to unify the treatment of both the Riemann zeta function and Dirichlet L-functions (for an in-depth exploration, consult Chapter 12 of \cite{Apostol}).

Appell polynomials \cite{Appell, Costabile}, which serve as a generalization of the classical Bernoulli, Euler, and Hermite polynomials, are formulated based on an arbitrary sequence of real numbers $a_{n}$. They can be expressed as $A_{n}(x)=\sum_{k=0}^{n}\binom{n}{k}a_{k} x^{n-k}$. In \cite{Navas.Appell, Navas.Anote}, L.M. Navas, F.J. Ruiz, and J.L. Varona extended the connection between Bernoulli polynomials and the Hurwitz zeta function to encompass a broader class of Appell polynomials and their associated functions. In this study, we explore the correspondence between these mathematical concepts from a distinctive and comprehensive standpoint. Specifically, our primary objective is to construct a category of functions, derived from Appell polynomials, which serve as a generalization of the Hurwitz zeta function and satisfy a functional equation analogous to the Hurwitz's formula presented below (see e.g., Apostol, chapter 12, page 256): \[\zeta(1-z,a)=\frac{\Gamma(z)}{(2\pi)^{z}}\left(e^{-\frac{i\pi z}{2}}F(a,z)+e^{\frac{i\pi z}{2}}F(-a,z)\right).\] where $F(s,z)$ is the periodic zeta function given by $F(s,z)=\sum_{n=1}^{+\infty}\frac{e^{2in\pi s}}{n^z}$, and $\Gamma(z)$ is the Euler gamma function. The latter formula represents a generalization of the functional equation for Riemann's zeta function. In the next presentation, we will disclose our key findings and clarify the underlying motivation for undertaking this investigation.

Consider a real-analytic function $f(t)=\sum_{n=0}^{+\infty}\frac{C_{f,n}}{n!}t^n$ at $0$ and set $p_{f}(t):=\frac{e^{t}-1}{t}f(t)$. The \textit{C-polynomials associated to} $f$ and the \textit{P-polynomials associated to} $f$ are defined by their respective exponential generating functions: \[f(t)e^{xt}=\sum_{n=0}^{+\infty}\frac{C_{f,n}(x)}{n!}t^n \ \text{ and }\ p_{f}(t)e^{xt}=\sum_{n=0}^{+\infty}\frac{P_{f,n}(x)}{n!}t^n.\] The two sequences of numbers, $C_{f,n}=C_{f,n}(0)$ and $P_{f,n}=P_{f,n}(0)$, are respectively referred to as \textit{C-numbers associated to} $f$ and \textit{P-numbers associated to} $f$. Consequently, we obtain: \[C_{f,n}(x)=\sum_{k=0}^{n} \binom{n}{k}C_{f,k} x^{n-k} \ \text{ and }\  P_{f,n}(x)=\sum_{k=0}^{n}\binom{n}{k}P_{f,k} x^{n-k}.\]

We observe that both C-polynomials and P-polynomials are Appell polynomials, defined by the sequences of coefficients $C_{f,n}$ and $P_{f,n}$, respectively. When we choose the function $f(t)$ as $\beta(t)=\frac{t}{e^{t}-1}$, the C-polynomials $C_{\beta,n}(x)$ reduce to the classical Bernoulli polynomials $B_{n}(x)$, which are generated by the function $\beta(t)e^{xt}$, and the P-polynomials $P_{\beta,n}(x)$ reduce to the monomials $x^n$, with $p_{\beta}(t)$ being the constant function equal to $1$. We note that the polynomial $P_{f,n}(x)$, denoted by $x^{(n,f)}$, generalizes the power $x^n$. This observation has motivated our exploration of a new generalization of the Hurwitz zeta function theory. Our investigation begins with a generalization of Faulhaber's formula (or sums of powers, see, e.g., \cite{Knuth}) for the finite sum over $j$ of the terms $j^{(n,f)}$, where $n\geq0$ is a fixed integer, and concludes with the introduction and analysis of the \textit{LC-function associated to} $f$, defined for $\Re(z)>1$ by the series: \[L(z,f):=\sum_{n=n_{f}}^{+\infty}n^{(-z,f)},\] where $s^{(z,f)}$ denotes the function of complex variables defined by: \[s^{(z,f)}:=\sum_{k=0}^{+\infty}\binom{z}{k}P_{f,k}s^{z-k}=s^{z}\sum_{k=0}^{+\infty}\binom{z}{k}P_{f,k}\left(\frac{1}{s}\right)^{k},\] and $n_{f}$ denotes the positive integer \[n_{f}:=\floor{\frac{1}{\rho_{f}}}+1.\] Here, $\rho_{f}$ is the convergence radius of the series $\sum_{k=0}^{+\infty}P_{f,k}s^{k}$, which we assume to be non-zero, and it is also the convergence radius of the series $\sum_{k=0}^{+\infty}\binom{z}{k}P_{f,k}s^{k}$ if $z\notin\mathbb{Z}_{\geq0}$. The function $s^{(z,f)}$ is a natural extension of the P-polynomials given above by $x^{(n,f)}=P_{f,n}(x)=\sum_{k=0}^{n}\binom{n}{k}P_{f,k} x^{n-k}$. For all integers $n\geq n_{f}$, since $\frac{1}{n}<\rho_{f}$, the function $z\mapsto n^{(z,f)}$ is well-defined on $\mathbb{C}$ and is even entire.
 
Note that when we choose the function $f(t)$ as $\beta_{a}(t)=te^{(a-1)t}/(e^{t}-1)$, where $0<a<1$, we obtain: \[p_{\beta_{a}}(t)=e^{(a-1)t}=\sum_{k=0}^{+\infty}\frac{(a-1)^{k}}{k!}t^{k}.\] Thus, the P-numbers associated to $\beta_{a}$ are $P_{\beta_{a},k}=(a-1)^k$ and $\rho_{\beta_{a}}=\frac{1}{1-a}>1$. As a result, $n_{\beta_{a}}=1$ and \[n^{(z,\beta_{a})}=n^{z}\sum_{k=0}^{+\infty}\binom{z}{k}\left(\frac{a-1}{n}\right)^{k}=(n+a-1)^{z}.\] Consequently, $L(z,\beta_{a})$ is none other than the Hurwitz zeta function $\zeta(z,a)$.\\

We now outline the structure of the paper and highlight the main results. In \Cref{sec2}, we discuss well-known properties of Appell polynomials and derive new ones, focusing on both C-polynomials and P-polynomials, as well as the identities connecting these two types of polynomials. We then present three formulas related to C-polynomials that generalize classical results, specifically the multiplication formula for Bernoulli polynomials (see e.g. \cite{Olver}, page 590), the Euler-Maclaurin formula (see e.g. \cite{Arakawa}, Chapter 5 and \cite{Olver}, page 63), and Faulhaber's formula. In \Cref{sec3}, we introduce and prove several useful lemmas, which will be applied in subsequent sections, particularly in \Cref{sec4}, to study the function $s^{(z,f)}$. \Cref{sec5}, the essential part of the article, is devoted to the LC-function $L(z,f)$ and its main properties. We first express this function using an integral representation on the interval $(0,+\infty)$, and then, by employing a second integral representation along a Hankel contour, we obtain the analytic continuation of $L(z,f)$ on the whole complex plane, except perhaps for a single simple pole at $1$, depending on the value of $f(0)=C_{f,0}$. At this stage, we can prove that the equation $L(-n,f)=-\frac{C_{f,n+1}(n_{f})}{n+1}$ holds for all integers $n\geq 0$. Finally, we derive a formula for the LC-function, which generalizes Hurwitz's formula.\\

\paragraph{\textbf{Notations}. } 
Throughout this paper, we will use the following notations:

\begin{enumerate}[1.]
\item \label{nota1} $\mathbb{Z}$, $\mathbb{Z}_{\geq1}$, and $\mathbb{Z}_{\geq0}$ represent, respectively, the set of integers, the set of positive integers, and the set of non-negative integers.
\item \label{nota2} The function $\beta(t)=\frac{t}{e^t-1}$ serves as the exponential generating function for the Bernoulli numbers $B_{n}$ of the first kind ($B_{1}=-\frac{1}{2}$),
\begin{equation}\label{eqt1.1}
\beta(t)=\sum_{n=0}^{+\infty} \frac{B_{n}}{n!}t^n, \quad |t|<2\pi.
\end{equation}
\item \label{nota3} The assignment symbol \qq{$:=$} is employed as a abbreviation for \qq{is defined as} or \qq{is denoted by}.
\item \label{nota4} If $f$ is a function, we use $\overline{f}$, $\underline{f}$, and $f_{(\alpha)}$ ($\alpha\in\mathbb{C}$) to represent the functions defined by $\overline{f}(t)=f(-t)$, $\underline{f}(t)=e^{-t}f(-t)$, and $f_{(\alpha)}(t)=\frac{e^{\alpha t}-1}{\alpha(e^{t}-1)}f(\alpha t)$, respectively.
\item \label{nota5} The notation $\floor{.}$ indicates the integer part function.
\end{enumerate}

\section{C-numbers and C-polynomials}\label{sec2}

During the whole paper, we assume that $f$ is a real-analytic function at $0$, represented by its Maclaurin series expansion, as given bellow: 
\begin{equation}\label{eqt2.1}
f(t)=\sum_{n=0}^{+\infty} \frac{C_{f,n}}{n!}t^n.
\end{equation}
In this section, our aim is to revisit and clarify various properties related to Appell polynomials, as discussed in \cite{Appell,Costabile}. Our investigation covers both well-known and newly discovered aspects of these polynomials.

\subsection{Definitions and main properties}\label{subsec2.1}

As observed in (\ref{eqt2.1}), $f$ is the exponential generating function of the coefficients $C_{f,n}$  which are referred to as \textit{C-numbers associated to} $f$. We introduce the \textit{C-polynomials associated to $f$}, denoted by $C_{f,n}(x)$, through the subsequent exponential generating function:
\begin{equation}\label{eqt2.2}
f(t)e^{xt}=\sum_{n=0}^{+\infty} \frac{C_{f,n}(x)}{n!}t^n.
\end{equation}
It is evident that for all $n\geq 0$,
\begin{equation}\label{eqt2.3}
C_{f,n}(0)=C_{f,n}.
\end{equation}
\begin{remark}\label{rem2.1}
If we choose the function $f(t)$ to be $\beta(t)=\frac{t}{e^{t}-1}$, we obtain the sequences of Bernoulli numbers and polynomials (see \hyperref[nota2]{Notations.2}). That is to say, $C_{\beta,n}(x)=B_{n}(x)$ and $C_{\beta,n}=B_{n}$.
\end{remark}

\begin{proposition}\label{propo2.1}
The C-polynomials satisfy, for every $n \geq 0$,
\begin{equation}\label{eqt2.4}
C_{f,n}(x)=\sum_{k=0}^{n} \binom{n}{k} C_{f,k} x^{n-k}.
\end{equation}
\end{proposition}
\begin{proof}
Utilizing the Cauchy product formula, we expand the function $t\mapsto f(t)e^{xt}$ in a power series about $0$ and substitute the obtained expression into (\ref{eqt2.2}) to derive,
\[
\sum_{n=0}^{+\infty} \frac{C_{f,n}(x)}{n!}t^n =\left( \sum_{k=0}^{+\infty} \frac{C_{f,k}}{k!}t^k \right) \left( \sum_{p=0}^{+\infty} \frac{x^p}{p!}t^p \right)=\sum_{n=0}^{+\infty}\left( \sum_{k=0}^{n} \binom{n}{k} C_{f,k} x^{n-k} \right) \frac{t^n}{n!}.\]
By equating the coefficients of the different powers of $t$ between the two series in the equality above, we achieve the desired result.
\end{proof}

\begin{example}\label{exmp2.1}
Here are the first three values of C-polynomials associated to $f$.
\begin{enumerate}[]
\item $C_{f,0}(x)=C_{f,0}$.
\item $C_{f,1}(x)=C_{f,0}x+C_{f,1}$.
\item $C_{f,2}(x)=C_{f,0}x^2+2C_{f,1}x+C_{f,2}$.
\end{enumerate}
\end{example}

\begin{proposition}\label{propo2.2}
The following properties hold for all C-polynomials associated to $f$:\\
$-$ For all $n \geq 1$,
\begin{equation}\label{eqt2.5}
C'_{f,n}(x)=nC_{f,n-1}(x),
\end{equation}
with $C'_{f,0}(x)=0$.\\
$-$ For all $n \geq 0$,
\begin{equation}\label{eqt2.6}
\int_x^{x+1} C_{f,n}(t) \diff{t} = \frac{C_{f,n+1}(x+1)-C_{f,n+1}(x)}{n+1}.
\end{equation}
\end{proposition}

\begin{proof}
(i) Let us differentiate both sides of Equation (\ref{eqt2.2}) with respect to $x$. This yields: \[\sum_{n=0}^{+\infty} \frac{C'_{f,n}(x)}{n!}t^n=tf(t)e^{xt}=\sum_{n=0}^{+\infty} \frac{C_{f,n}(x)}{n!}t^{n+1}=\sum_{n=1}^{+\infty} \frac{nC_{f,n-1}(x)}{n!}t^n.\] By matching the coefficients on both sides of the equation, we obtain the desired result.

(ii) From (\ref{eqt2.5}) it follows that for all $n \geq 0$, \[\int_x^{x+1} C_{f,n}(t)\diff{t} = \frac{1}{n+1} \int_x^{x+1} C'_{f,n+1}(t) \diff{t} = \frac{C_{f,n+1}(x+1)-C_{f,n+1}(x)}{n+1}.\]
\end{proof}

\begin{proposition}\label{propo2.3}
Consider the two functions $\overline{f}(t) = f(-t)$ and $\underline{f}(t) = e^{-t}f(-t)$, as defined in \hyperref[nota4]{Notations.$4$}. For all non-negative integers $n$, the two following equations hold:
\begin{equation}\label{eqt2.7}
C_{f,n}(-x)=(-1)^nC_{\overline{f},n}(x).
\end{equation}
\begin{equation}\label{eqt2.8}
C_{f,n}(1-x)=(-1)^n C_{\underline{f},n}(x).
\end{equation}
\end{proposition}

\begin{proof}
Using (\ref{eqt2.4}), the first equation can be obtained as follows: \[C_{f,n}(-x)=(-1)^{n}\sum_{k=0}^{n}\binom{n}{k}(-1)^kC_{f,k}x^{n-k}=(-1)^{n}\sum_{k=0}^{n}\binom{n}{k}C_{\overline{f},k}x^{n-k}.\] Based on (\ref{eqt2.2}), we have: \[\sum_{n=0}^{+\infty} \frac{C_{f,n}(1-x)}{n!}t^n = f(t)e^{(1-x)t} = \underline{f}(-t)e^{-tx} = \sum_{n=0}^{+\infty} \frac{C_{\underline{f},n}(x)}{n!}(-t)^n.\]
By comparing the coefficients on both sides, we derive the second equation.
\end{proof}

\begin{proposition}\label{propo2.4}
For all non-negative integers $n$ and for all real numbers $x$ and $y$, we have the following equations:
\begin{equation}\label{eqt2.9}
C_{f,n}(x+y)=\sum_{k=0}^{n}\binom{n}{k}C_{f,k}(x)y^{n-k}.
\end{equation}
\begin{equation}\label{eqt2.10}
C_{f,n}(x+1)=\sum_{k=0}^{n}\binom{n}{k}C_{f,k}(x).
\end{equation}
\end{proposition}

\begin{proof}
In order to prove the desired result, we employ (\ref{eqt2.2}) in conjunction with the Cauchy product formula, as demonstrated below:
\begin{align*}
\sum_{n=0}^{+\infty}\frac{C_{f,n}(x+y)}{n!}t^n&=\left(\sum_{k=0}^{+\infty}\frac{C_{f,k}(x)}{k!}t^k\right)\left(\sum_{p=0}^{+\infty}\frac{(ty)^p}{p!}\right)\\
&=\sum_{n=0}^{+\infty}\left(\sum_{k=0}^{n}\binom{n}{k}C_{f,k}(x)y^{n-k}\right)\frac{t^n}{n!}.
\end{align*}
\end{proof}

\subsection{P-numbers and P-polynomials}\label{subsec2.2}
In this subsection, we present a special instance of C-polynomials, which subsequently enables the generalization of monomials $x^n$ at the end of subsection \ref{subsec2.3}.

We consider the function
\begin{equation} \label{eqt2.11}
p_f(t) := \frac{e^t - 1}{t}f(t),
\end{equation}
which is real-analytic at $0$. We denote the \textit{P-polynomials associated to} $f$ as $P_{f,n}(x)$, defined using the exponential generating function:
\begin{equation}\label{eqt2.12}
p_f(t)e^{xt} = \sum_{n=0}^{+\infty}\frac{P_{f,n}(x)}{n!}t^n.
\end{equation}
We define the corresponding \textit{P-numbers associated to} $f$ as
\begin{equation}\label{eqt2.13}
P_{f,n} = P_{f,n}(0),
\end{equation}
which, by taking $x=0$ in (\ref{eqt2.12}), gives
\begin{equation}\label{eqt2.14}
p_f(t) = \sum_{n=0}^{+\infty}\frac{P_{f,n}}{n!}t^n.
\end{equation}

By comparing Equations (\ref{eqt2.12}) and (\ref{eqt2.2}), we observe that the P-numbers and the P-polynomials associated to the function $f$ are, respectively, the C-numbers and the C-polynomials associated to the function $p_f$. In other words, $P_{f,n}(x) = C_{p_f,n}(x)$ and $P_{f,n} = C_{p_f,n}$. We can therefore easily deduce the properties of the P-polynomials based on what has already been proved in Subsection \ref{subsec2.1} for C-polynomials. Consequently, for all $n \geq 0$, and for all $x,y \in \mathbb{R}$, the following formulas hold:
\begin{equation}\label{eqt2.15}
P_{f,n}(x) = \sum_{k=0}^{n}\binom{n}{k}P_{f,k}x^{n-k},
\end{equation}
with $P_{f,0}(x)=P_{f,0}$.
\begin{equation}\label{eqt2.16}
P'_{f,n}(x) = nP_{f,n-1}(x),  \quad (n\geq0)
\end{equation}
with $P'_{f,0}(x) = 0$.
\begin{equation}\label{eqt2.17}
\int_x^{x+1} P_{f,n}(t) \diff{t} = \frac{P_{f,n+1}(x+1) - P_{f,n+1}(x)}{n+1}.
\end{equation}
\begin{equation}\label{eqt2.18}
P_{f,n}(x+y) = \sum_{k=0}^{n}\binom{n}{k}P_{f,k}(x)y^{n-k}.
\end{equation}

\subsection{Some relations between C-polynomials and P-polynomials} \label{subsec2.3}

The two propositions presented below provide a set of formulas that express, respectively, the P-polynomials and P-numbers in terms of their counterparts, the C-polynomials and C-numbers, and vice versa. These equations will facilitate a deeper understanding and practical application of these concepts, ultimately enabling the derivation of significant results in subsequent analyses.

\begin{proposition}\label{propo2.5}
For all non-negative integers $n$, the P-polynomials can be expressed in terms of C-polynomials according to the following expression:
\begin{equation}\label{eqt2.19}
P_{f,n}(x)=\frac{1}{n+1}\sum_{k=0}^{n} \binom{n+1}{k}C_{f,k}(x).
\end{equation}
Furthermore, when $x=0$, the expression simplifies to:
\begin{equation}\label{eqt2.20}
P_{f,n}=\frac{1}{n+1}\sum_{k=0}^{n}\binom{n+1}{k}C_{f,k}.
\end{equation}
\end{proposition}

\begin{proof}
By considering Equations (\ref{eqt2.12}), (\ref{eqt2.11}), (\ref{eqt2.2}), and the Cauchy product formula, we obtain,
\begin{align*}
\sum_{n=0}^{+\infty}\frac{P_{f,n}(x)}{n!}t^n = p_f(t)e^{xt}&=\frac{1}{t}e^t\left(f(t)e^{xt}\right)-\frac{1}{t}f(t)e^{xt} \\ 
&=\frac{1}{t}\left(\sum_{k=0}^{+\infty}\frac{t^k}{k!}\right)\left(\sum_{p=0}^{+\infty}\frac{C_{f,p}(x)}{p!}t^p\right) - \frac{1}{t}f(t)e^{xt}\\
&=\frac{1}{t}\sum_{n=0}^{+\infty}\left(\sum_{k=0}^{n} \binom{n}{k}C_{f,k}(x)\right)\frac{t^n}{n!} - \frac{1}{t}\left(\sum_{n=0}^{+\infty}\frac{C_{f,n}(x)}{n!}t^n\right)\\
&=\frac{1}{t}\sum_{n=0}^{+\infty}\left(\sum_{k=0}^{n} \binom{n}{k}C_{f,k}(x) - C_{f,n}(x)\right)\frac{t^n}{n!} \\
&=\frac{1}{t}\sum_{n=1}^{+\infty}\left(\sum_{k=0}^{n-1} \binom{n}{k}C_{f,k}(x)\right)\frac{t^n}{n!} \\
&=\sum_{n=0}^{+\infty}\left(\frac{1}{n+1}\sum_{k=0}^{n} \binom{n+1}{k}C_{f,k}(x)\right)\frac{t^n}{n!}.
\end{align*}
This leads to the desired result.
\end{proof}

\begin{samepage}
\begin{proposition}\label{propo2.6}
For all non-negative integers $n$, the C-polynomials can be expressed in terms of both, the P-polynomials and the Bernoulli numbers, according to the following expression:
\begin{equation} \label{eqt2.21}
C_{f,n}(x)=\sum_{k=0}^{n} \binom{n}{k}B_{n-k}P_{f,k}(x)=\sum_{k=0}^{n}\binom{n}{k}B_{k}P_{f,n-k}(x),
\end{equation}
Furthermore, evaluating the expression at $x=0$, we obtain
\begin{equation} \label{eqt2.22}
C_{f,n}=\sum_{k=0}^{n}\binom{n}{k}B_{n-k}P_{f,k}=\sum_{k=0}^{n}\binom{n}{k}B_{k}P_{f,n-k}.
\end{equation}
\end{proposition}
\end{samepage}

\begin{proof}
Referring to (\ref{eqt2.11}), we have: \[f(t)e^{xt}=\frac{t}{e^t-1}p_f(t)e^{xt}.\]
Consequently, by using Equations (\ref{eqt2.2}), (\ref{eqt1.1}), (\ref{eqt2.12}), and employing the Cauchy product formula, we obtain the following expression:
\begin{align*}
\sum_{n=0}^{+\infty}\frac{C_{f,n}(x)}{n!}t^n&=\left(\sum_{k=0}^{+\infty}\frac{B_k}{k!}t^k\right)\left(\sum_{p=0}^{+\infty}\frac{P_{f,n}(x)}{p!}t^p\right) \\
&=\sum_{n=0}^{+\infty}\left(\sum_{k=0}^{n}\binom{n}{k}B_{n-k}P_{f,k}(x)\right)\frac{t^n}{n!}.
\end{align*}
This leads to the formula (\ref{eqt2.21}).
\end{proof}

The main idea behind this paper manifests in the fundamental identity (\ref{eqt2.24}) pertaining to the C-polynomials, the proof of which is provided in the subsequent \Cref{propo2.7}. This equation offers a generalization of the next well-established fundamental identity associated with the classical Bernoulli polynomials:
\begin{equation}\label{eqt2.23}
B_{n}(x+1)-B_{n}(x)=nx^{n-1},\quad (n\geq 1).
\end{equation}

\begin{proposition}\label{propo2.7}
The fundamental identity of the C-polynomials can be expressed for all $n \geq 1$, as follows:
\begin{equation} \label{eqt2.24}
C_{f,n}(x+1)-C_{f,n}(x)=nP_{f,n-1}(x),
\end{equation}
\end{proposition}
\begin{proof}
Evidently, based on Equations (\ref{eqt2.10}) and (\ref{eqt2.19}), we can deduce that\[C_{f,n}(x+1)-C_{f,n}(x)=\sum_{k=0}^{n-1}\binom{n}{k}C_{f,k}(x)=nP_{f,n-1}(x).\]
\end{proof}

\begin{remark}\label{rem2.2}
\begin{enumerate}
\item Utilizing (\ref{eqt2.6}), we remark that the fundamental identity expressed in (\ref{eqt2.24}) can be represented in an alternative form for all $n \geq 0$ as:
\begin{equation}\label{eqt2.25}
\int_x^{x+1}C_{f,n}(t)\diff{t}=P_{f,n}(x).
\end{equation}
\item By setting $x=0$ in (\ref{eqt2.24}), we obtain the following formula for all $n \geq 1$:
\begin{equation}\label{eqt2.26}
C_{f,n}(1)=C_{f,n}+nP_{f,n-1}.
\end{equation}
\end{enumerate}
\end{remark}

By considering the function $\beta(t)=\frac{t}{e^{t}-1}$ as $f$ in (\ref{eqt2.24}), and in accordance with \Cref{rem2.1} and (\ref{eqt2.23}), we obtain for all $n\geq 1$: \[B_{n}(x+1)-B_{n}(x)=nP_{\beta,n-1}(x)=nx^{n-1}.\] This result implies that $P_{\beta,n}(x)=x^{n}$. A comparison of Equations (\ref{eqt2.23}) and (\ref{eqt2.24}) reveals that the P-polynomials $P_{f,n}(x)$ generalize the simple powers $x^n$. This remark justifies the choice of the letter \qq{P} in the notation of these polynomials. In order to become familiar with this generalization, we prefer to represent these polynomials using the following notation:
\begin{equation}\label{eqt2.27}
x^{(n,f)}:=P_{f,n}(x),
\end{equation}
where $x^{(n,\beta)}=x^n$. Motivated by this observation, we propose that the theory of the Hurwitz zeta function, initially developed through the study of Bernoulli polynomials, can be generalized by examining the properties of the C-polynomials $C_{f,n}(x)$ (which generalize Bernoulli polynomials) and the P-polynomials $x^{(n,f)}$ (which generalize the monomials $x^n$).

The first step involves establishing a formula for calculating the finite sum of the terms $j^{(n,f)}$ over $j$ (refer to Subsection \ref{subsec2.4}). However, a more challenging task is to define and study a function of complex variables $P_{f}(s,z)$, such that for $(z,s):=(x,n)\in\mathbb{R}\times\mathbb{Z}_{\geq0}$, it coincides with the P-polynomial $x^{(n,f)}$. In other words, $P_{f}(x,n)=P_{f,n}(x)=x^{(n,f)}$. In \Cref{sec3}, we will denote this function as $s^{(z,f)}$ and use it in \Cref{sec4} to generalize the Hurwitz zeta function.

We conclude this \Cref{sec2} by generalizing the three prominent formulas: namely, Faulhaber's formula, the multiplication formulas (also referred to as the distribution formula), and, finally, the Euler-Maclaurin formula.

\subsection{Generalized Faulhaber formula}\label{subsec2.4}

Faulhaber's formula \cite{Knuth} provides a method to calculate the sum of the $n$-th powers of the first $m$ non-negative integers, with $n \geq 0$. This formula includes the Bernoulli numbers of the first kind ($B_{1} = -\frac{1}{2}$) and can be expressed as: \[\sum_{j=0}^{m-1}j^{n}=\frac{1}{n+1}\sum_{k=0}^{n}\binom{n+1}{k}B_{k}m^{n+1-k}.\] In this subsection, we aim to expand upon this classical formula by introducing a more comprehensive version. The proposed formula calculates the finite sum $\sum_{j=0}^{m-1}j^{(n,f)}$ and incorporates the C-numbers associated to $f$.

\begin{theorem}
For all $m \in \mathbb{Z}_{\geq 1}$ and for all $n \in \mathbb{Z}_{\geq 0}$, We introduce two general forms for the proposed formula:

\begin{enumerate}[$-$]
\item First form:
\begin{equation}\label{eqt2.28}
\sum_{j=0}^{m-1}(x+j)^{(n,f)}=\frac{C_{f,n+1}(x+m)-C_{f,n+1}(x)}{n+1}.
\end{equation}

\item Second form:
\begin{equation}\label{eqt2.29}
\sum_{j=0}^{m-1}(x+j)^{(n,f)}=\frac{1}{n+1}\sum_{k=0}^{n}\binom{n+1}{k}C_{f,k}(x)m^{n+1-k}.
\end{equation}
\end{enumerate}
\end{theorem}

\begin{proof}
Using (\ref{eqt2.24}) and applying a telescopic sum, we obtain:
\begin{align*}
\sum\limits_{j=0}^{m-1}(x+j)^{(n,f)}&=\sum\limits_{j=0}^{m-1}{{P_{f,n}}(x+j)}\\
&=\sum\limits_{j=0}^{m-1}{\frac{{{C_{f,n+1}}(x+j+1)-{C_{f,n+1}}(x+j)}}{{n+1}}}\\
&=\frac{C_{f,n+1}(x+m)-C_{f,n+1}(x)}{n+1}.
\end{align*}

Furthermore, from (\ref{eqt2.9}), we have:
\begin{align*}
C_{f,n+1}(x+m)&=\sum_{k=0}^{n+1}\binom{n+1}{k}C_{f,k}(x)m^{n+1-k}\\
&=\sum_{k=0}^{n}\binom{n+1}{k} C_{f,k}(x)m^{n+1-k}+C_{f,n+1}(x).
\end{align*}

Thus, we obtain the anticipated outcome.
\end{proof}

Upon setting $x=0$ in Equation (\ref{eqt2.29}), we derive the generalized Faulhaber formula, proposed in the next corollary:
 
\begin{corollary}[Generalized Faulhaber Formula]
For all $m \in \mathbb{Z}_{\geq 1}$ and for all $n \in \mathbb{Z}_{\geq 0}$,
\begin{equation}\label{eqt2.30}
\sum_{j=0}^{m-1}j^{(n,f)}=\frac{1}{n+1}\sum_{k=0}^{n}\binom{n+1}{k} C_{f,k} m^{n+1-k}.\\
\end{equation}
\end{corollary}

\subsection{Multiplication formula for the C-polynomials}\label{subsec2.5}

The Raabe's multiplication formula for the Bernoulli polynomials, originally presented by Joseph Ludwig Raabe \cite{Raabe1851} in 1851, is given by: \[\sum_{k=0}^{m-1}B_{n}\left(\frac{x+k}{m}\right)=\frac{1}{m^{n-1}}B_{n}(x).\] It should be noted that not all Appell sequences necessarily adhere to this specific formula. A more general result was provided by Carlitz in his research \cite{Carlitz1953TheMF}. In the subsequent theorem, we propose an alternative approach for deriving a more comprehensive multiplication formula suitable for all Appell polynomials $C_{f,n}(x)$.

\begin{theorem}[Multiplication formula for the C-polynomials]\label{thm2.3}
For all integers $n \geq 0$ and $m \geq 1$, we can express the following relation:
\begin{equation}\label{eqt2.31}
\sum_{k=0}^{m-1}C_{f,n}\left(\frac{x + k}{m} \right)=\frac{1}{m^{n-1}}C_{f_{(m)},n}(x),
\end{equation}
where the function $f_{(m)}(t)$ is defined as: $\displaystyle f_{(m)}(t)=\frac{e^{mt}-1}{m(e^t-1)}f(mt)$. Based on $(\ref{eqt2.11})$, it can be readily verified that: $\displaystyle f_{(m)}(t)=\beta(t)p_f(mt)$ $($see \hyperref[nota2]{Notations.$2$}$)$.
\end{theorem}

\begin{proof}
We define the following expression:
\[M_{f,n,m}:=m^{n}\sum_{k=0}^{m - 1}C_{f,n}\left(\frac{x+k}{m}\right).\] The proof involves combining (\ref{eqt2.2}) with the formula for the sum of a geometric series. The following series expansion can be expressed as:
\begin{align*}
\sum_{n=0}^{\infty}\frac{M_{f,n,m}}{n!} t^n &= \sum_{k=0}^{m-1}\sum_{n=0}^{\infty} \frac{C_{f,n}\left(\frac{x+k}{m}\right)}{n!}(mt)^n 
= \sum_{k=0}^{m-1}f(mt)e^{(x+k)t} \\
&= f(mt)e^{xt}\frac{e^{mt}-1}{e^{t}-1} = mf_{(m)}(t)e^{xt} = m\sum_{n=0}^{\infty}\frac{C_{f_{(m)},n}(x)}{n!} t^n.
\end{align*}
To complete the proof, we must compare the coefficients of the two series.
\end{proof}

\subsection{Euler-Maclaurin formula in terms of C-numbers and P-numbers} \label{subsec2.6}

In this subsection, we aim to rewrite the Euler-Maclaurin formula \cite{Arakawa,Olver} using C-numbers and P-numbers in place of Bernoulli numbers.

Consider $n \in \mathbb{Z}_{\geq 0}$. The function $x \mapsto C_{f,n}(x - \floor{x})$, defined for all $x \in \mathbb{R}$, where $\floor{x}$ denotes the integer part of $x$, fulfills the condition ${C}_{f,n}\left(x-\floor{x}\right)=C_{f,n}\left(x-k\right)$, for any integer $k$ and any real number $x$ such that $k \leq x < k + 1$.

\begin{theorem}\label{thm.}
Let $f(0) \neq 0$, which, according to $(\ref{eqt2.1})$, is equivalent to $C_{f,0} \neq 0$. Let $m$ and $n$ be integers, and consider a function $g$ that is $N$ times continuously differentiable on the interval $[m, n]$, where $N \geq 2$. Then, the following equation holds:
\begin{align*}
\sum_{k=m}^n g(k)&=\int_m^{n}g(x)\diff{x}+\frac{C_{f,1}(1)g(m)-C_{f,1}g(n)}{C_{f,0}}\\
&\quad+\frac{1}{C_{f,0}}\sum_{r=2}^N\frac{(-1)^r C_{f,r}}{r!}\left(g^{(r-1)}(n)-g^{(r-1)}(m)\right)\\
&\quad +\frac{1}{C_{f,0}}\sum_{r=2}^N \frac{(-1)^rP_{f,r-1}}{(r-1)!}\sum_{k=m+1}^{n}g^{(r-1)}(k)+R_{N},
\end{align*}
where \[R_N=\frac{(-1)^{N+1}}{C_{f,0}N!} \int_m^{n} {C}_{f,N}\left(x-\floor{x}\right)g^{(N)}(x)\diff{x}.\] In the formula above, according to $(\ref{eqt2.26})$ and $(\ref{eqt2.20})$, we can express $C_{f,1}(1)$ in terms of the C-numbers $C_{f,0}$ and $C_{f,1}$ as follows: $C_{f,1}(1) = C_{f,1} + P_{f,0} = C_{f,1} + C_{f,0}$.
\end{theorem}

\begin{proof}
Set for $r \geq 0$ and $k\in\mathbb{Z}$: \[I_{r,k}:=\int_{k}^{k+1}C_{f,r}(x-\floor{x})g^{(r)}(x)\diff{x}.\]
Apply the change of variables $t=x-k$, followed by integration by parts to obtain:
\begin{align*}
I_{r+1,k}&=\int_{k}^{k+1}C_{f,r+1}(x-k)g^{(r+1)}(x)\diff{x}=\int_{0}^{1}C_{f,r+1}(t)g^{(r+1)}(t+k)\diff{x}\\
&= \left[C_{f,r+1}(t)g^{(r)}(t+k)\right]_{0}^{1}-(r+1)I_{r,k}.
\end{align*}
Summing over the range $k = m, \dots, n-1$ and multiply both sides by $\frac{(-1)^r}{(r+1)!C_{f,0}}$,
\begin{align*}
R_{r+1}-R_{r}&=\frac{(-1)^r}{(r+1)!C_{f,0}}\left[C_{f,r+1}(t)\sum_{k=m}^{n-1}g^{(r)}(t+k)\right]_{0}^{1}.
\end{align*}
Summing over $r = 0, \dots, N-1$, we obtain: \[R_{N}-R_{0}=\frac{1}{C_{f,0}}\sum_{r=0}^{N-1}\frac{(-1)^r}{(r+1)!} \left(C_{f,r+1}(1)\sum_{k=m}^{n-1} g^{(r)}(k+1) - C_{f,r+1}\sum_{k=m}^{n-1} g^{(r)}(k) \right),\]
with $R_0=-\int_{m}^{n}g(x)\diff{x}$. From (\ref{eqt2.26}), we have $C_{f,r+1}(1)=C_{f,r+1}+(r+1)P_{f,r}$. Thus, we can reorganize the terms as follows:
\begin{align*}
R_{N}-R_{0}&=\frac{1}{C_{f,0}}\sum_{r=0}^{N-1}\frac{(-1)^rC_{f,r+1}}{(r+1)!}\left(\sum_{k=m+1}^{n} g^{(r)}(k)-\sum_{k=m}^{n-1}g^{(r)}(k)\right)\\
&\quad +\frac{1}{C_{f,0}}\sum_{r=0}^{N-1}\frac{(-1)^rP_{f,r}}{r!}\sum_{k=m+1}^{n}g^{(r)}(k)\\
&=\frac{1}{C_{f,0}}\sum_{r=0}^{N-1}\frac{(-1)^rC_{f,r+1}}{(r+1)!}\left(g^{(r)}(n)-g^{(r)}(m)\right)\\
&\quad +\frac{1}{C_{f,0}}\sum_{r=0}^{N-1}\frac{(-1)^rP_{f,r}}{r!}\sum_{k=m+1}^{n}g^{(r)}(k)
\end{align*}
Using (\ref{eqt2.20}) to substitute $P_{f,0}$ with $C_{f,0}$, we can further simplify the expression:
\begin{align*}
R_{N}-R_{0}&=\frac{C_{f,1}\left(g(n)-g(m)\right)}{C_{f,0}}+\frac{1}{C_{f,0}}\sum_{r=1}^{N-1}\frac{(-1)^rC_{f,r+1}}{(r+1)!}\left(g^{(r)}(n)-g^{(r)}(m)\right)\\
&\quad +\sum_{k=m+1}^{n}g(k)+\frac{1}{C_{f,0}}\sum_{r=1}^{N-1}\frac{(-1)^rP_{f,r}}{r!}\sum_{k=m+1}^{n}g^{(r)}(k)\\
&=\sum_{k=m}^{n}g(k)+\frac{C_{f,1}g(n)-(C_{f,0}+C_{f,1})g(m)}{C_{f,0}}\\
&\quad + \frac{1}{C_{f,0}}\sum_{r=1}^{N-1}\frac{(-1)^rC_{f,r+1}}{(r+1)!}\left(g^{(r)}(n)-g^{(r)}(m)\right)\\
&\quad +\frac{1}{C_{f,0}}\sum_{r=1}^{N-1}\frac{(-1)^rP_{f,r}}{r!}\sum_{k=m+1}^{n}g^{(r)}(k).
\end{align*}
This concludes the proof.
\end{proof}

\section{Some useful lemmas}\label{sec3}

In this section, we present three fundamental lemmas that will be utilized in subsequent discussions.

Let $(u_n)_{n\geq0}$ be a sequence of real numbers such that the power series $\sum_{n=0}u_{n}s^{n}$ has a non-zero radius of convergence $\rho$, which may potentially be infinite.

\begin{samepage}
\begin{lemma}\label{lem3.1}
Let us fix $z\in\mathbb{C}$.
\begin{enumerate}
\item
\begin{enumerate}[$(a)$]
\item If $z$ is a non negative integer, then the power series $\sum_{n=0}\binom{z}{n}u_{n}s^{n}$ converges everywhere.
\item If $z\in\mathbb{C}\backslash\mathbb{Z}_{\geq0}$, then the radius of convergence of the power series $\sum_{n=0}\binom{z}{n}u_{n}s^{n}$ is $\rho$.
\end{enumerate}
\item  The power series $\sum_{n=0}\frac{u_n}{n!}t^n$ has an infinite radius of convergence and its sum defines an infinitely differentiable function $u$ on $\mathbb{R}$. According to $(\ref{eqt2.1})$, the numbers $u_n$ are the C-numbers $C_{u,n}$ associated to $u$.
\end{enumerate}
\end{lemma}
\end{samepage}

\begin{proof}
\begin{enumerate}
\item We will prove the result by considering two distinct cases: when $z \in \mathbb{Z}_{\geq 0}$ and when $z \notin \mathbb{Z}_{\geq 0}$.
\begin{enumerate}
\item If $z\in\mathbb{Z}_{\geq0}$, all terms $\binom{z}{n}$ become zero starting from rank $z+1$. Thus, for all $s\in\mathbb{C}$, we have \[\sum_{n=0}^{+\infty}\binom{z}{n}u_{n}s^{n}=\sum_{n=0}^{z}\binom{z}{n}u_{n}s^{n}.\]
\item We now assume that $z\notin\mathbb{Z}_{\geq0}$. For all $n\geq 1$, we have
\begin{equation}\label{eqt3.1}
\left\lvert\binom{z}{n} u_n\right\rvert^{\frac{1}{n}}= \left\lvert\binom{z}{n}\right\rvert^{\frac{1}{n}}|u_n|^{\frac{1}{n}}.
\end{equation}
Applying Cesaro's theorem, we get \[\frac{1}{n}\ln\left\lvert\binom{z}{n}\right\rvert = \frac{1}{n}\sum_{n=0}^{n-1}\ln\left\lvert\frac{z-k}{k+1}\right\rvert\xrightarrow[\ n \ ]{}0,\] and thus, by continuity of the exponential function at $0$, $\left\lvert\binom{z}{n}\right\rvert^{\frac{1}{n}}\xrightarrow[ n ]{}1$. The equality (\ref{eqt3.1}), due to the positivity of terms, leads to \[\limsup_{n}\left\lvert\binom{z}{n} u_n\right\rvert^{\frac{1}{n}} = \limsup_{n}|u_n|^{\frac{1}{n}}.\] Therefore, by the Cauchy-Hadamard theorem, the two series $\sum_{n=0}\binom{z}{n}u_{n}s^{n}$ and $\sum_{n=0}u_{n}s^{n}$ have the same radius of convergence $\rho$.
\end{enumerate}
\item It is evident that
$\lim_{n}\left(\frac{1}{n!}\right)^{\frac{1}{n}}=0$ and $\limsup_{n}|u_n|^{\frac{1}{n}}<+\infty$. Thus, \[\limsup_{n}\left\lvert\frac{u_n}{n!}\right\rvert^{\frac{1}{n}}=\limsup_{n}\left(\frac{1}{n!}\right)^{\frac{1}{n}}|u_n|^{\frac{1}{n}}=0.\]
Consequently, by the Cauchy-Hadamard theorem, we obtain the desired result.
\end{enumerate}
\end{proof}

\begin{figure}[H]
\begin{center}
\includegraphics[]{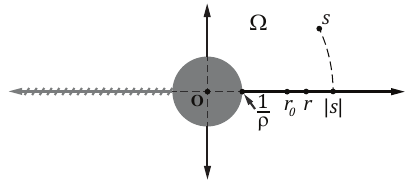}
\caption{Figure illustrating the definition domain $\Omega$ for the function $s\mapsto S_u(s,z)$, with fixed $z$. The domain $\Omega$ is composed of the complex plane, omitting the negative real axis and the shaded gray disk.}
\label{fig:Fig1}
\end{center}
\end{figure}

We define $\Omega$ as follows (\Cref{fig:Fig1}):
\begin{equation}\label{eqt3.2}
\Omega := \left\{ s \in \mathbb{C}\backslash\mathbb{R}_{\leq 0},\ |s| > \frac{1}{\rho}\right\}.
\end{equation}
Now, according to \Cref{lem3.1}.(i), we consider the function defined for every $s \in \Omega$ and every $z \in \mathbb{C}$ by
\begin{equation}\label{eqt3.3}
S_u(s,z):=\sum_{n=0}^{+\infty}\binom{z}{n}u_ns^{z-n}=s^{z}\sum_{n=0}^{+\infty}\binom{z}{n}u_n\left(\frac{1}{s}\right)^{n},
\end{equation}
where $u$ is the function mentioned in \Cref{lem3.1}.(ii).

\begin{lemma}\label{lem3.2}
\begin{enumerate}
\item For a fixed $s \in \Omega$, $S_u(s, z)$ is an entire function of $z$.
\item Let $s \in \Omega$ and $z \in \mathbb{C}$. Let $r$ and $r_0$ be arbitrary real numbers such that $\frac{1}{\rho} < r_0 < r \leq |s|$ $($\Cref{fig:Fig1}$)$. Then,
\begin{enumerate}[$(a)$]
\item If $|z| \leq M$ $(M>0)$, there exists a constant $K_{M}$ depending on $M$ and independent of $z$, such that
\begin{equation}\label{eqt3.5}
\left\lvert S_u(s, z)\right\rvert \leq \frac{rK_{M}}{r-r_0}\left\lvert s^{z}\right\rvert.
\end{equation}
\item There exists a constant $K_{z}$ depending on $z$, such that
\begin{equation}\label{eqt3.4}
\left\lvert S_u(s, z)\right\rvert \leq \frac{rK_{z}}{r-r_0}\left\lvert s^{z}\right\rvert.
\end{equation}
\end{enumerate}
\end{enumerate}
\end{lemma}

\begin{proof}
\begin{enumerate}
\item Let us fix $s \in \Omega$. We first note that the function of $z$,
\begin{equation}\label{eqt3.6}
S_u(s,z) = s^z \sum_{n=0}^{+\infty} \binom{z}{n} u_{n}s^{-n},
\end{equation}
is the product of an entire function and the sum of a series of entire functions. We consider an arbitrary compact disk $|z| \leq M$ and prove that this series converges uniformly on every such disk. This will establish that the sum $\sum_{n=0}^{+\infty} \binom{z}{n} u_{n}s^{-n}$ defines an entire function of $z$, and therefore the function $z \mapsto S_u(s, z)$ is entire.

From the inequality
\begin{equation}\label{eqt3.7}
\left\lvert\binom{z}{n}\right\rvert \leq \left\lvert\binom{-M}{n}\right\rvert,
\end{equation}
which can be verified using the triangle inequality and the two identities \[\binom{z}{n}=\frac{\prod_{k=0}^{n-1}(z-k)}{n!}\text{ and }\binom{-z}{n}=\frac{(-1)^n\prod_{k=0}^{n-1}(z+k)}{n!},\] we obtain, for all integers $n \geq 0$, \[\left\lvert\binom{z}{n}u_{n}s^{-n}\right\rvert \leq \left\lvert\binom{-M}{n}u_{n}\left(\frac{1}{r_0}\right)^{n}\right\rvert\left(\frac{r_0}{|s|}\right)^n.\] By \Cref{lem3.1}.(i), since $\frac{1}{r_0} < \rho$, the sequence $\left(\binom{-M}{n}u_{n}\left(\frac{1}{r_0}\right)^{n}\right)_{n}$ is bounded by a constant $K_M > 0$ that depends on $M$ but not on $z$. Thus, for all integers $n \geq 0$,
\begin{equation}\label{eqt3.8}
\left\lvert\binom{z}{n}u_{n}s^{-n}\right\rvert\leq K_M\left(\frac{r_0}{|s|}\right)^n.
\end{equation}
As $0 < \frac{r_0}{|s|} < 1$, we conclude that the series $\sum_{n=0} \binom{z}{n} u_{n}s^{-n}$ converges uniformly on every compact disk $|z| \leq M$. Hence, the result follows.
\item 
\begin{enumerate}[(a)]
\item From (\ref{eqt3.6}) and (\ref{eqt3.8}), we derive \[\left\lvert S_u(s,z)\right\rvert\leq\frac{K_{M}}{1-\frac{r_0}{|s|}}\left\lvert s^{z}\right\rvert.\] Consequently, since $\frac{r_0}{|s|} \leq \frac{r_0}{r} < 1$, \[\left\lvert S_u(s,z)\right\rvert\leq\frac{K_{M}}{1-\frac{r_0}{r}}\left\lvert s^{z}\right\rvert=\frac{rK_{M}}{r-r_0}\left\lvert s^{z}\right\rvert.\]
\item Similarly to the approach taken above, we have for all integers $n \geq 0$, \[\left\lvert\binom{z}{n}u_{n}s^{-n}\right\rvert=\left\lvert\binom{z}{n}u_{n}\left(\frac{1}{r_0}\right)^{n}\right\rvert\left(\frac{r_0}{|s|}\right)^n.\] According to \Cref{lem3.1}.(i), since $\frac{1}{r_0} < \rho$, the sequence $\left(\binom{z}{n}u_{n}\left(\frac{1}{r_0}\right)^{n}\right)_n$ is bounded by a constant $K_z > 0$, dependent on $z$. Thus, for all integers $n \geq 0$, \[\left\lvert\binom{z}{n}u_{n}s^{-n}\right\rvert\leq K_{z}\left(\frac{r_0}{|s|}\right)^n.\] As a result, since $\frac{r_0}{|s|} \leq \frac{r_0}{r} < 1$, \[\left\lvert S_u(s,z)\right\rvert\leq\frac{K_{z}}{1-\frac{r_0}{|s|}}\left\lvert s^{z}\right\rvert\leq\frac{K_{z}}{1-\frac{r_0}{r}}\left\lvert s^{z}\right\rvert.\]
\end{enumerate}
\end{enumerate}
This concludes the proof.
\end{proof}

The subsequent lemma enables us to express the function $S_u(s, z)$ in an integral form, utilizing two well-known formulas of the gamma function:
\begin{equation}\label{eqt3.9}
\binom{-z}{n}=\frac{(-1)^n\Gamma(z+n)}{\Gamma(z)n!}, \ z\in\mathbb{C}\backslash\mathbb{Z}_{\leq0},\ n\in\mathbb{Z}_{\geq0}
\end{equation}
and
\begin{equation}\label{eqt3.10}
s^{-z}\Gamma(z)=\int_0^{+\infty}t^{z-1}e^{-st}\diff{t}, \ \Re(z)>0,\ \Re(s)>0.
\end{equation}
For all the subsequent analysis, we set $z = \sigma + iy$, where $\sigma$ and $y$ denote the real and imaginary parts of $z$, respectively. It is worth noting that the forthcoming lemma represents a particular instance of Ramanujan's master theorem \cite{Amdeberhan}.

\begin{samepage}
\begin{lemma}\label{lem3.3}
Let $\sigma > 0$ and $\Re(s) > \frac{1}{\rho}$. Then the following statements hold: 
\begin{enumerate}
\item The integral $\int_0^{+\infty}t^{z-1}e^{-st}u(-t)\diff{t}$ converges and 
\begin{equation}\label{eqt3.11}
S_u(s,-z)=\frac{1}{\Gamma(z)}\int_0^{+\infty}t^{z-1}e^{-st}u(-t)\diff{t},
\end{equation}
which can be expressed as
\begin{equation}\label{eqt3.12}
\int_0^{+\infty}t^{z-1}e^{-st}\left(\sum_{n=0}^{+\infty}\frac{u_{n}}{n!}(-t)^n\right)\diff{t}=\Gamma(z)\sum_{n=0}^{+\infty}\binom{-z}{n}u_{n}s^{-z-n}.
\end{equation}

\item In $(i)$, let us consider the function $\overline{u}(t):=u(-t)$ as the choice for $u(t)$. Given that the convergence radius of the power series $\sum_{n=0}\overline{u}_{n}s^{n}$ is also $\rho$, the integral $\int_0^{+\infty}t^{z-1}e^{-st}u(t)\diff{t}$ converges. Consequently, we have:
\begin{equation}\label{eqt3.16}
S_{\overline{u}}(s,-z)=\frac{1}{\Gamma(z)}\int_{0}^{+\infty}t^{z-1}e^{-st}u(t)\diff{t}.
\end{equation}
This means that,
\begin{equation}\label{eqt3.12}
\int_0^{+\infty}t^{z-1}e^{-st}\left(\sum_{n=0}^{+\infty}\frac{u_{n}}{n!}t^{n}\right)\diff{t}=\Gamma(z)\sum_{n=0}^{+\infty}\binom{-z}{n}(-1)^{n}u_{n}s^{-z-n}.
\end{equation}
\end{enumerate}
\end{lemma}
\end{samepage}

\begin{proof} To start, let us fix $z, s \in \mathbb{C}$ with the conditions $\sigma > 0$ and $\Re(s) > \frac{1}{\rho}$. Using Equations (\ref{eqt3.3}), (\ref{eqt3.9}), and (\ref{eqt3.10}), we can transform $S_u(s, -z)$ as shown below:
\begin{align*}
S_u(s,-z)&=\frac{1}{\Gamma(z)}\sum_{n=0}^{\infty}\frac{(-1)^nu_{n}s^{-(z+n)}\Gamma(z+n)}{n!}\\
&=\frac{1}{\Gamma(z)}\sum_{n=0}^{\infty}\int_0^{+\infty}\frac{t^{z+n-1}e^{-st}(-1)^nu_{n}}{n!}\diff{t}.
\end{align*}
Thus, we write:
\begin{equation}\label{eqt3.13}
S_u(s,-z)=\frac{1}{\Gamma(z)}\sum_{n=0}^{\infty} \int_0^{+\infty}g_n(t)\diff{t},
\end{equation}
where $g_n$ represents the functions of $t$:
\begin{equation}\label{eqt3.14}
g_n(t)=\frac{t^{z+n-1}e^{-st}(-1)^nu_{n}}{n!}.
\end{equation}
Next, we need to interchange the sum and integral in (\ref{eqt3.13}). The simplest way to achieve this is to perform the term-by-term integration. Functions $g_n$ are continuous on the interval $(0,+\infty)$, and as stated in \Cref{lem3.1}.(ii), \[\sum_{n=0}^{+\infty}g_n(t) = t^{z-1}e^{-st}\sum_{n=0}^{+\infty}\frac{u_{n}}{n!}(-t)^n=t^{z-1}e^{-st}u(-t).\] Consequently, the series of functions $\sum_{n=0}g_n(t)$ converges on $(0,+\infty)$ to the function $t \mapsto t^{z-1}e^{-st}u(-t)$, which is continuous on $(0,+\infty)$ according to \Cref{lem3.1}.(ii). Now, it remains to verify the convergence of the series $\sum_{n=0}\int_{0}^{+\infty}\left\lvert g_n(t)\right\rvert\diff{t}$. Indeed, \[\int_{0}^{+\infty}\left\lvert g_n(t)\right\rvert\diff{t}=\frac{|u_{n}|}{n!}\int_{0}^{+\infty}t^{\sigma+n-1}e^{-\Re(s) t}\diff{t},\] which we can rewrite, considering (\ref{eqt3.10}) and (\ref{eqt3.9}), since $\sigma+n>0$ and $\Re(s)>0$, \[\int_{0}^{+\infty}\left\lvert g_n(t)\right\rvert\diff{t}=\frac{|u_{n}|\Gamma(\sigma+n)}{n!\Re(s)^{\sigma+n}}
=\frac{\Gamma(\sigma)}{\Re(s)^{\sigma}}\binom{-\sigma}{n}|u_{n}|\left(-\frac{1}{\Re(s)}\right)^{n}.\] According to \Cref{lem3.1}.(i), since $0<\frac{1}{\Re(s)}<\rho$, the series $\sum_{n=0}\int_{0}^{+\infty}\left\lvert g_n(t)\right\rvert\diff{t}$ converges. As a result, by applying the term-by-term integration theorem, the function $t \mapsto t^{z-1}e^{-st}u(-t)$ is integrable on $(0,+\infty)$, and (\ref{eqt3.13}) becomes \[S_{u}(s,-z)=\frac{1}{\Gamma(z)}\int_0^{+\infty}\sum_{n=0}^{\infty}g_n(t)\diff{t}=\frac{1}{\Gamma(z)}\int_0^{+\infty}t^{z-1}e^{-st}u(-t)\diff{t}.\] This completes the proof.
\end{proof}

\section{The function $s^{(z,f)}$}\label{sec4}

From now on, we presume that the real-analytic function $f(t) = \sum_{n=0}^{\infty} \frac{C_{f,n}}{n!} t^n$ is centered at $0$ such that the series $\sum_{n=0}P_{f,n}s^{n}$ has a non-zero radius of convergence $\rho_{f}\leq+\infty$, where $P_{f,n}$ are the P-numbers associated to $f$ as introduced in Subsection \ref{subsec2.2}.

Upon substituting the sequence $u_n$ with the sequence $P_{f,n}$ in Lemma \ref{lem3.1}, the following conclusions can be deduced.

\begin{theorem}\label{thm4.1}
Fix $z\in\mathbb{C}$.
\begin{enumerate}
\item If $z\in\mathbb{C}\backslash\mathbb{Z}_{\geq0}$, then the radius of convergence of the power series $ \sum_{n=0}\binom{z}{n}P_{f,n}s^{n}$ is given by $\rho_{f}$.
\item The power series
$\sum_{n=0}\frac{P_{f,n}}{n!}t^n$ has an infinite radius of convergence and its sum $p_{f}(t)=\frac{e^{t}-1}{t}f(t)$, seen in $(\ref{eqt2.14})$ and $(\ref{eqt2.11})$, is a function of class $\mathcal{C}^{\infty}$ on $\mathbb{R}$.
\end{enumerate}
\end{theorem}

\begin{figure}[H]
\begin{center}
\includegraphics[]{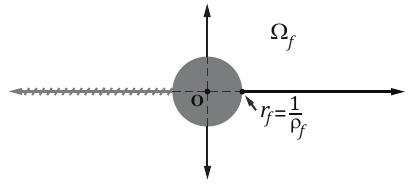}
\caption{Figure depicting the domain of definition $\Omega_{f}$ for the function $s\mapsto P_{f}(s,z)$, where $z$ is held constant. The domain $\Omega_{f}$ consists of the complex plane, excluding the negative real axis and the shaded gray disk.}
\label{fig:Fig2}
\end{center}
\end{figure}

Let $r_{f}:=\frac{1}{\rho_{f}}$ with $r_{f}=0$ if $\rho_{f}=+\infty$. As we have done in the previous \Cref{sec3}, we define $\Omega_{f}$ as shown in \Cref{fig:Fig2}:
\begin{equation}\label{eqt4.1}
\Omega_{f}:=\left\{s\in \mathbb{C}\mid s\notin\mathbb{R}_{\leq 0},\ |s|>r_{f}\right\},
\end{equation}
We then consider a function, $P_{f}(s,z)$, defined for all $s\in\Omega_{f}$ and $z\in\mathbb{C}$ as follows:
\begin{equation}\label{eqt4.2}
P_{f}(s,z):=S_{p_{f}}(s,z).
\end{equation}
Using (\ref{eqt3.3}), we can express the function as:
\begin{equation}\label{eqt4.3}
P_{f}(s,z)=\sum_{n=0}^{+\infty}\binom{z}{n}P_{f,n}s^{z-n}=s^{z}\sum_{n=0}^{+\infty}\binom{z}{n}P_{f,n}\left(\frac{1}{s}\right)^{n}.
\end{equation}

Let $N\in \mathbb{Z}_{\geq0}$. Observe that the function of $s$, denoted as \[P_{f}(s,N):=\sum_{n=0}^{N}\binom{N}{n}P_{f,n}s^{N-n},\] coincides with the P-polynomial $P_{f,N}(s)$ associated to $f$, as outlined in Subsection \ref{subsec2.2} and specified by (\ref{eqt2.15}). Utilizing the notation in (\ref{eqt2.27}), we write $s^{(N,f)}=P_{f,N}(s)=P_{f}(s,N)$. Consequently, we can extend this notation to the function $P_{f}(s,z)$, as illustrated below.
\begin{equation}\label{eqt4.4}
s^{(z,f)}:=P_{f}(s,z)=\sum_{n=0}^{+\infty}\binom{z}{n}P_{f,n}s^{z-n},\quad s\in\Omega_{f},\ z\in\mathbb{C}.
\end{equation}

According to (\ref{eqt4.2}) and the established notation above, we can write $s^{(z,f)}=S_{p_{f}}(s,z)$. Consequently, the two outcomes presented in the subsequent theorem emerge as direct corollaries of \Cref{lem3.2}.(i) and \Cref{lem3.3}.(i).

\begin{theorem}\label{thm4.2}
\begin{enumerate}
\item For every $s\in\Omega_{f}$, the function $s^{(z,f)}$ is entire with respect to $z$.
\item For any $\sigma>0$ and $\Re(s)>r_{f}$, the integral $\int_0^{+\infty}t^{z-1}e^{-st}p_{f}(-t)\diff{t}$ converges, and we have
\begin{equation}\label{eqt4.5}
s^{(-z,f)}=\frac{1}{\Gamma(z)}\int_0^{+\infty}t^{z-1}e^{-st}p_{f}(-t)\diff{t},
\end{equation}
which implies
\begin{equation}\label{eqt4.6}
\int_0^{+\infty}t^{z-1}e^{-st}\left(\sum_{n=0}^{+\infty}\frac{P_{f,n}}{n!}(-t)^n\right)\diff{t}=\frac{\Gamma(z)}{s^{z}}\sum_{n=0}^{+\infty}\binom{-z}{n}P_{f,n}s^{-n}.\\
\end{equation}
\end{enumerate}
\end{theorem}

For the remainder of this discussion, we introduce the following notations and properties. Let $\alpha\in\mathbb{C}\backslash\{0\}$, and consider the function
\begin{equation}\label{eqt4.7}
f_{(\alpha)}(t):=\frac{e^{\alpha t}-1}{\alpha(e^{t}-1)}f(\alpha t).
\end{equation}
We observe that
\begin{equation}\label{eqt4.8}
f_{(1)}(t)=f(t) \text{ and } \underline{f}(t):=f_{(-1)}(t)=e^{-t}f(t),
\end{equation}
which implies
\begin{equation}\label{eqt4.9}
\underline{\underline{f}}(t)=f(t).
\end{equation}
By (\ref{eqt2.11}),
\begin{equation}\label{eqt4.10}
p_{f_{(\alpha)}}(t)=p_{f}(\alpha t) \text{ and } p_{\underline{f}}(t)=p_{f_{(-1)}}(t)=p_{f}(-t).
\end{equation}
Consequently, by (\ref{eqt2.14}), the first equation in (\ref{eqt4.10}) becomes \[\sum_{n=0}^{+\infty}\frac{P_{f_{(\alpha)},n}}{n!}t^n=\sum_{n=0}^{+\infty}\frac{\alpha^{n}P_{f,n}}{n!}t^n.\] Comparing coefficients of both sums, we deduce for all $n\geq0$,
\begin{equation}\label{eqt4.11}
P_{f_{(\alpha)},n}=\alpha^{n}P_{f,n} \text{ and } P_{\underline{f},n}=P_{f_{(-1)},n}=(-1)^n P_{f,n}.
\end{equation}
Since $\sum_{n=0}P_{f_{(\alpha)},n}s^{n}=\sum_{n=0}P_{f,n}(\alpha s)^{n}$, the convergence radius $\rho_{f_{(\alpha)}}$ of the series on the left-hand side is given by
\begin{equation}\label{eqt4.12}
\rho_{f_{(\alpha)}}=\frac{\rho_{f}}{|\alpha|} \text{ with } \rho_{\underline{f}}=\rho_{f_{(-1)}}=\rho_{f}.
\end{equation}
Thus, we define
\begin{equation}\label{eqt4.13}
r_{f_{(\alpha)}}:=\frac{1}{\rho_{f_{(\alpha)}}}=|\alpha|r_{f} \text{ with } r_{\underline{f}}=r_{f}.
\end{equation}
By (\ref{eqt4.4}) and (\ref{eqt4.11}), for all $s\in\Omega_{f_{(\alpha)}}$ and for all $z\in\mathbb{C}$,
\begin{equation}\label{eqt4.14}
s^{(z,f_{(\alpha)})}=\sum_{n=0}^{+\infty}\binom{z}{n}\alpha^{n}P_{f,n}s^{z-n}.
\end{equation}
This yields, for all $\beta\in\mathbb{C}\backslash\{0\}$, 
\begin{equation}\label{eqt4.15}
s^{(z,f_{(\alpha)})}=\beta^{z}\sum_{n=0}^{+\infty}\binom{z}{n}\left(\frac{\alpha}{\beta}\right)^{n} P_{f,n}\left(\frac{s}{\beta}\right)^{z-n}=\beta^{z}\left(\frac{s}{\beta}\right)^{(z,f_{(\alpha/\beta)})}.
\end{equation}
\newline 

As an immediate consequence of \Cref{thm4.2}.(ii), the following corollary is obtained.

\begin{corollary}\label{cor4.1}
\begin{enumerate}
\item For all $\sigma>0$ and for all $\Re(s)>r_{f_{(\alpha)}}$,
\begin{equation}\label{eqt4.16}
s^{(-z,f_{(\alpha)})}=\frac{1}{\Gamma(z)}\int_0^{+\infty}t^{z-1}e^{-st}p_{f_{(-\alpha)}}(t)\diff{t},
\end{equation}
\item Taking $\alpha=-1$ in $\text{(i)}$, we obtain for all $\sigma>0$ and for all $\Re(s)>r_{f}$,
\begin{equation}\label{eqt4.17}
s^{(-z,\underline{f})}=\frac{1}{\Gamma(z)}\int_0^{+\infty}t^{z-1}e^{-st}p_{f}(t)\diff{t},
\end{equation}
where, according to $(\ref{eqt4.8})$, $\underline{f}=f_{(-1)}$,  $f_{(1)}=f$ and $r_{f_{(-1)}}=r_{f}$.
\end{enumerate}
\end{corollary} 

\section{The LC-Function} \label{sec5}

\subsection{Introduction}\label{subsec5.1}

In this section, the majority of proofs are influenced by those presented in (\cite{Apostol}, Chapter 12, pages 251-259) concerning the Hurwitz zeta function, albeit with some technical differences.\

\begin{figure}[H]
\begin{center}
\includegraphics[]{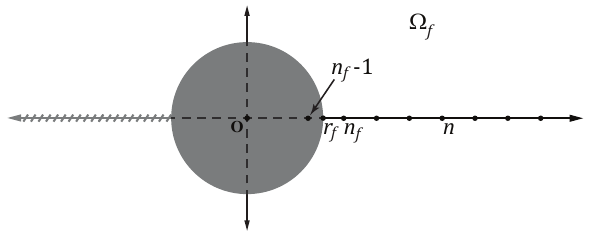}
\caption{$n_{f}$ is the smallest integer in the domain $\Omega_{f}$, defined by $n_{f} := \left\lfloor{r_{f}}\right\rfloor + 1$.}
\label{fig:Fig3}
\end{center}
\end{figure}

Consider $n_{f}$ to be the smallest integer within the set $\Omega_{f}$ (\Cref{fig:Fig3}). Based on (\ref{eqt4.1}), we can define $n_{f}$ as follows:
\begin{equation}\label{eqt5.1}
n_{f}:=\left\lfloor{r_{f}}\right\rfloor+1=\left\lfloor{\frac{1}{\rho_{f}}}\right\rfloor+1.
\end{equation}

\begin{remark}\label{rem5.1}
If we take an integer $n\geq n_{f}$, we can see from (\ref{eqt5.1}) that $n> r_{f}$. This indicates that $n$ is an element of $\Omega_{f}$ (\Cref{fig:Fig3}), as per (\ref{eqt4.1}). Therefore, it can be inferred from \Cref{thm4.2}.(i) that the function $z\mapsto n^{(-z,f)}$ is entire.
\end{remark}

In this section we introduce and investigate the \textit{LC-function associated to the function} $f$, defined by the series 
\begin{equation}\label{eqt5.2}
L(z,f):=\sum_{n=n_{f}}^{+\infty}n^{(-z,f)}.
\end{equation}

\begin{theorem}\label{thm5.1}
The series defining the function $L(z,f)$ converges absolutely for $\sigma>1$. The convergence is uniform within every compact subset of the half-plane $\sigma>1$. Furthermore, given that $n^{(-z,f)}$ $(n\geq n_{f})$ is a sequence of entire functions of $z$, as demonstrated by \Cref{rem5.1}, $L(z,f)$ is an analytic function of $z$ within the half-plane $\sigma>1$.
\end{theorem}

\begin{proof}
Let $r_{0}$ be an arbitrary real number with $r_{f}\leq r_{0}<n_{f}$ and $K$ be a compact subset of the half-plane $\sigma>1$. Then, there exist two constants $M>1$ and $\delta>1$ such that, for all $z\in K$, $|z|<M$ and $\sigma\geq \delta$. From (\ref{eqt4.4}) and (\ref{eqt4.2}) we have, $n^{(-z,f)}=S_{p_{f}}(n,-z)$. According to \Cref{lem3.2}.(ii).(a), since $r_{f}<r_{0}<n_{f}\leq n$ ($r_f=\frac{1}{\rho_{f}}$), there exists a positive constant $K_{M}$, dependent on $M$ and not on $z$, such that \[\left\lvert n^{(-z,f)}\right\rvert\leq\frac{n_{f}K_{M}}{n_{f}-r_{0}}n^{-\sigma}.\]
Since $\sigma>1$, all the statements of the theorem follows from the inequalities 
\[\sum_{n=n_{f}}^{+\infty}\left\lvert n^{(-z,f)}\right\rvert \leq \frac{n_{f}K_{M}}{n_{f}-r_{0}}\sum_{n=n_{f}}^{+\infty}\frac{1}{n^{\sigma}}\leq \frac{n_{f}K_{M}}{n_{f}-r_{0}}\sum_{n=n_{f}}^{+\infty}\frac{1}{n^{\delta}}.\]
\end{proof}

\subsection{Integral representation for the LC-function}\label{subsec5.2}

\begin{theorem}\label{thm5.2}
For $\sigma>1$, the integral representation is given by
\begin{equation}\label{eqt5.3}
L(z,f)=\frac{1}{\Gamma(z)}\int_0^{+\infty}\frac{t^{z-1}e^{(1-n_{f})t}p_{f}(-t)}{e^t-1}\diff{t}=\frac{1}{\Gamma(z)}\int_0^{+\infty}\frac{t^{z-1}e^{(1-n_{f})t}p_{\underline{f}}(t)}{e^t-1}\diff{t},
\end{equation}
which we can rewrite, using $(\ref{eqt2.11})$, as 
\begin{equation}\label{eqt5.4}
L(z,f)=\frac{1}{\Gamma(z)}\int_0^{+\infty}t^{z-2}e^{(1-n_{f})t}\underline{f}(t)\diff{t}.
\end{equation}
\end{theorem}

\begin{remark}\label{rem5.2}
The function $\Gamma(z)L(z,f)$ represents the Mellin transform of the function $\varphi_{f}(t)=t^{-1}e^{(1-n_{f})t}\underline{f}(t)$. Thus, we can express the LC-function as
\begin{equation}\label{eqt5.5}
L(z,f)=\frac{1}{\Gamma(z)}{\mathcal{M}\varphi_{f}}(z).
\end{equation}
\end{remark}

\begin{proof}[Proof of $\Cref{thm5.2}$]
Let $\sigma>1$ and $n\geq n_{f}$ ($n_{f}>r_{f}$). According to \Cref{thm4.2}.(ii), we have \[\Gamma(z)n^{(-z,f)}=\int_0^{+\infty}t^{z-1}e^{-nt}p_{f}(-t)\diff{t}.\]
Now, summing over the integers $n\geq n_{f}$ and using (\ref{eqt5.2}), it comes
\begin{equation}\label{eqt5.6}
\Gamma(z)L(z,f)=\sum_{n=n_{f}}^{+\infty}\int_{0}^{+\infty}t^{z-1}e^{-nt}p_{f}(-t)\diff{t}=\sum_{n=n_{f}}^{+\infty}\int_{0}^{+\infty}h_{n}(t)\diff{t},
\end{equation}
where, for every $n\geq n_{f}$, $h_n(t):=t^{z-1}e^{-nt}p_{f}(-t)$. In order to interchange the sum and the integral, we apply the term-by-term integration theorem. First, according to \Cref{thm4.1}.(ii), $h_n(t)$ forms a sequence of continuous functions on $(0,+\infty)$. Since
\begin{equation}\label{eqt5.7}
\sum_{n=n_{f}}^{+\infty}h_n(t)=\frac{t^{z-1}e^{-n_{f}t}p_{f}(-t)}{1-e^{-t}},
\end{equation}

the series of functions $\sum_{n=n_{f}} h_n(t)$ converges on $(0,+\infty)$ to a continuous function. Now, it remains to verify the convergence of the series $\sum_{n=n_{f}}\int_{0}^{+\infty}\left\lvert h_n(t)\right\rvert\diff{t}$. By using (\ref{eqt2.14}), we write
 \[\int_{0}^{+\infty}\left\lvert h_n(t)\right\rvert\diff{t}=\int_{0}^{+\infty}t^{\sigma-1} e^{-nt}\left\lvert p_{f}(-t)\right\rvert\diff{t}\leq\int_{0}^{+\infty}t^{\sigma-1}e^{-nt}\sum_{k=0}^{+\infty}\frac{|P_{f,k}|}{k!}t^k\diff{t}.\]

From \Cref{lem3.3}.(ii), since $\sigma>1$ and $n>r_{f}=\frac{1}{\rho_{f}}$ ($\rho_f$ is also the radius of convergence of the series $\sum_{k=0}\left\lvert P_{f,k}\right\rvert s^{k}$), the integral on the right-hand side of the previous inequality converges, with 
\[\int_{0}^{+\infty}t^{\sigma-1}e^{-nt}\sum_{k=0}^{+\infty}\frac{|P_{f,k}|}{k!}t^k \diff{t}=\frac{\Gamma(\sigma)}{n^{\sigma}} \sum_{k=0}^{\infty}\binom{-\sigma}{k}(-1)^k|P_{f,k}|n^{-k}.\] 
Since $\binom{-\sigma}{k}(-1)^k = \frac{\Gamma(\sigma+k)}{\Gamma(\sigma)k!}>0$, 
\[\int_{0}^{+\infty}\left\lvert h_n(t)\right\rvert\diff{t}\leq\frac{\Gamma(\sigma)}{n^{\sigma}}\sum_{k=0}^{\infty}\left\lvert\binom{-\sigma}{k} P_{f,k}\right\rvert n^{-k}.\]
Let us choose an arbitrary real number $r$ such that $\frac{1}{\rho_f}<r<n_{f}$. By \Cref{thm4.1}.(i), the sequence $\left( \left\lvert\binom{-\sigma}{k}P_{f,k}r^{-k}\right\rvert\right)_{k}$ is therefore bounded by a constant $Q_{\sigma}$. Consequently, 
\[\int_{0}^{+\infty}\left\lvert h_n(t)\right\rvert\diff{t}\leq \frac{\Gamma(\sigma)}{n^{\sigma}} \sum_{k=0}^{\infty}\left\lvert\binom{-\sigma}{k} P_{f,k} \left(\frac{1}{r}\right)^k \right\rvert \left(\frac{r}{n}\right)^k \leq \frac{Q_{\sigma}\Gamma(\sigma)}{n^{\sigma}}\sum_{k=0}^{\infty}\left(\frac{r}{n}\right)^{k}.\]
Since $r<n_{f}\leq n$, we obtain \[\int_{0}^{+\infty}\left\lvert h_n(t)\right\rvert\diff{t}\leq\frac{Q_{\sigma}\Gamma(\sigma)}{n^{\sigma}}\frac{1}{1-\frac{r}{n}}\leq\frac{Q_{\sigma}\Gamma(\sigma)}{n^{\sigma}}\frac{1}{1-\frac{r}{n_{f}}}.\] Thus, \[\sum_{n=n_{f}}^{+\infty}\int_{0}^{+\infty}\left\lvert h_n(t)\right\rvert\diff{t}\leq\frac{n_{f}Q_{\sigma}\Gamma(\sigma)}{n_{f}-r}\zeta(\sigma,n_{f}).\] The result is obtained from (\ref{eqt5.7}) and (\ref{eqt4.10}), after performing term-to-term integration in (\ref{eqt5.6}), as follows: \[\Gamma(z)L(z,f)=\int_0^{+\infty }\sum_{n=n_{f}}^{+\infty} h_n(t)\diff{t}=\int_0^{+\infty}\frac{t^{z-1}e^{-n_{f}t}p_{\underline{f}}(t)}{1-e^{-t}}\diff{t}.\]
\end{proof}

\subsection{A contour integral representation for the LC-function}\label{subsec5.3}

\begin{figure}[H]
\begin{center}
\includegraphics[]{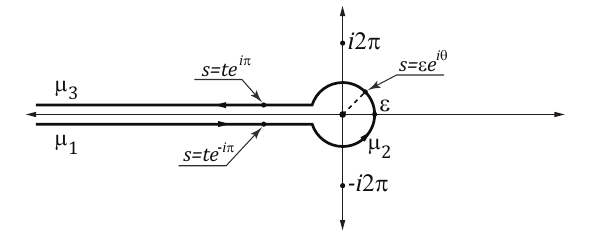}
\caption{The Hankel contour path $\mu$, traversed in the positive sense, with $\varepsilon < 2\pi$.}
\label{fig:Fig4}
\end{center}
\end{figure}

In order to extend the LC-function $L(z,f)$ to the half-plane $\sigma < 1$, an alternative integral representation is introduced, which is predicated on a positively oriented contour $\mu$. This contour consists of three distinct parts, denoted as $\mu_{1}$, $\mu_{2}$, and $\mu_{3}$, which are illustrated in \Cref{fig:Fig4}, with $\varepsilon < 2\pi$.

\begin{theorem}\label{thm5.3}
The function defined by the contour integral
\[I(z,f)=\frac{1}{2i\pi}\int_{\mu}\frac{s^{z-1}e^{n_{f}s}p_{f}(s)}{1-e^{s}}\diff{s}\]
is an entire function of $z$ and we have for $\sigma>1$ and $z\not=2,3,4...$,
\begin{equation}\label{eqt5.8}
L(z,f)=\Gamma(1-z)I(z,f).
\end{equation}
\end{theorem}

\begin{proof}
We commence by observing that on $\mu_{1}$ and $\mu_{3}$, we have $s^{z-1}=-t^{z-1}e^{-i\pi z}$ and $s^{z-1}=-t^{z-1}e^{i\pi z}$, respectively. Consequently, we can express $I(z,f)$ as
\begin{align*}
2i\pi I(z,f)&=\left(\int_{\mu_{1}}+\int_{\mu_{2}}+\int_{\mu_{3}}\right)\frac{s^{z-1}e^{n_{f}s}p_{f}(s)}{1-e^{s}}\diff{s}\\
&=e^{-i\pi z}\int_{+\infty}^{\varepsilon}\frac{t^{z-1}e^{-n_{f}t}p_{f}(-t)}{1-e^{-t}}\diff{t}+\int_{\mu_{2}}\frac{s^{z-1}e^{n_{f}s}p_{f}(s)}{1-e^{s}}\diff{s}\\
&\quad +e^{i\pi z}\int_{\varepsilon}^{+\infty}\frac{t^{z-1}e^{-n_{f}t}p_{f}(-t)}{1-e^{-t}}\diff{t}.
\end{align*}
From (\ref{eqt4.10}), which states that $p_{f}(-t)=p_{\underline{f}}(t)$, we can rewrite the above expression as
\begin{align}
2i\pi I(z,f)&=2i\sin(\pi z)\int_{\varepsilon}^{+\infty}\frac{t^{z-1}e^{(1-n_{f})t}p_{\underline{f}}(t)}{e^{t}-1}\diff{t}+\int_{\mu_{2}}\frac{s^{z-1}e^{n_{f}s}p_{f}(s)}{1-e^{s}}\diff{s}\label{eqt5.9}\\ 
&=2i\sin(\pi z)\int_{\varepsilon}^{+\infty}v(z,t)\diff{t}+I_{1}(z,\varepsilon).\label{eqt5.10}\\ 
&=2i\sin(\pi z)\int_{1}^{+\infty}v(z,t)\diff{t}+2i\sin(\pi z)\int_{\varepsilon}^{1}v(z,t)\diff{t}+I_{1}(z,\varepsilon).\nonumber
\end{align}
To demonstrate that $I(z,f)$ is an entire function, it suffices to apply the theorem of holomorphy under the integral sign to the integral $\int_{1}^{+\infty}v(z,t)\diff{t}$. Since $v(z,t)$, which is a continuous function of $t$ on $(1,+\infty)$, is an entire function of $z$, it remains to show that $v(z,t)$ is dominated on any compact disk $|z|\leq M$, $M>0$, by an integrable function of $t$ not depending on $z$. 

For $t\geq1$, we have \[|v(z,t)|=\frac{t^{\sigma-1}e^{(1-n_{f})t}}{e^{t}-1}\left\lvert p_{\underline{f}}(t)\right\rvert\leq\frac{e}{e-1}\left(t^{M-1}e^{-n_{f}t}\sum_{n=0}^{+\infty}\frac{\left\lvert P_{f,n}\right\rvert}{n!}t^n\right).\]
From (\ref{eqt5.1}), we have $n_{f}>\frac{1}{\rho_{f}}$. Hence, according to \Cref{lem3.3}.(ii), the function of $t$ in the right-hand side of the above inequality is integrable on $[1,+\infty)$.\\

Now, we proceed to prove (\ref{eqt5.8}). From (\ref{eqt5.10}), we have \[I_{1}(z,\varepsilon)=i\varepsilon^{z-1}\int_{-\pi}^{\pi}e^{i\theta(z-1)}v_{1}\left(\varepsilon e^{i\theta}\right)\diff{\theta},\] where \[v_{1}(s):=\frac{se^{n_{f}s}p_{f}(s)}{1-e^{s}}.\] Let $d$ be a fixed real number such that $\varepsilon<d<2\pi$. The function $v_{1}(s)$ being continuous on the compact disk $|s|\leq d$, it is bounded there by a constant A dependent only on $d$. Therefore, \[\left\lvert I_{1}(z,\varepsilon)\right\rvert\leq A\varepsilon^{\sigma-1}\int_{-\pi}^{\pi}e^{-y\theta}\diff{\theta},\] where $y$ is the imaginary part of $z$. Thus, for $\sigma>1$, we obtain \[\lim_{\varepsilon\to 0}I_{1}(z,\varepsilon)=0.\] By passing to the limit when $\varepsilon$ tends to $0$ in (\ref{eqt5.10}) and according to the integral representation of $L(z,f)$ in (\ref{eqt5.3}), we find \[2i\pi I(z,f)=2i\sin(\pi z)\int_{0}^{+\infty}\frac{t^{z-1}e^{(1-n_{f})t}p_{\underline{f}}(t)}{e^{t}-1}\diff{t}=2i\sin(\pi z)\Gamma(z)L(z,f).\] Hence, (\ref{eqt5.8}) follows from the complement formula $\Gamma(z)\Gamma(1-z)=\pi/\sin(\pi z)$.
\end{proof}

\subsection{The analytic continuation of the LC-function}\label{subsec5.4}

\begin{theorem}\label{thm5.3}
Let $f(0)\neq0$, then the function $L(z, f)$ has an analytic continuation to $\mathbb{C}\backslash\{1\}$ with a simple pole at $1$ and residue $f(0)$. On the other hand, if $f(0)=0$, $L(z, f)$ is an entire function. We have for all $z\not=1$,
\begin{equation}\label{eqt5.11}
L(z, f) = \Gamma(1 - z)I(z, f) = \frac{\Gamma(1 - z)}{2i\pi}\int_{\mu}\frac{s^{z - 1}e^{n_{f}s}p_{f}(s)}{1 - e^{s}}\diff{s}.
\end{equation}
\end{theorem}

\begin{proof}
Recall Formula (\ref{eqt5.8}), which is valid for $\sigma > 1$: \[L(z,f)=\Gamma(1-z)I(z,f).\] According to \Cref{thm5.3}, the function $I(z,f)$ is entire, so the function $\Gamma(1 - z)I(z, f)$ is analytic on the whole complex plane except possibly at the poles of $\Gamma(1 - z)$, namely, the integers $n \geq 1$. However, from \Cref{thm5.1}, the function $L(z, f)$ is analytic at all integers $n \geq 2$. Thus, the function $\Gamma(1 - z)I(z, f)$ is analytic on $\mathbb{C}\backslash\{1\}$. Consequently, the equation above provides the analytic continuation of $L(z, f)$ to the entire complex plane, except perhaps for $z = 1$.

Now, we proceed to examine the behavior of the function at the point $z = 1$. According to (\ref{eqt5.9}), we have
\begin{align}
I(1,f)&=\frac{1}{2i\pi}\int_{\mu_{2}}\frac{e^{n_{f}s}p_{f}(s)}{1-e^{s}}\diff{s}=\underset{s=0}{\text{Res}}\frac{e^{n_{f}s}p_{f}(s)}{1-e^{s}}\nonumber\\
&=\lim_{s\to 0}\frac{se^{n_{f}s}p_{f}(s)}{1-e^{s}}=-p_{f}(0)=-f(0).\label{eqt5.12}
\end{align}
$-$ If $f(0)\neq0$, the function $\Gamma(1-z)I(z,f)$ admits a simple pole at $1$, with residue \[I(1,f)\Res_{z=1}\Gamma(1-z)=f(0).\]
$-$ If $f(0)=0$, then by (\ref{eqt5.12}), $I(1,f)=0$. Since $I(z,f)$ is analytic at $1$ and $I(1, f) = 0$, the function $I(z, f) / (z - 1)$ is also analytic at $1$. Therefore, the function $L(z, f)$, which can be written as \[L(z,f)=\Gamma(1-z)I(z,f)=-(1-z)\Gamma(1-z)\frac{I(z,f)}{z-1}=-\Gamma(2-z)\frac{I(z,f)}{z-1},\] is analytic at $1$, with \[L(1,f)=-I'(1,f).\]
This completes the proof.
\end{proof}

\begin{remark}\label{rem5.3}
In the preceding proof, we demonstrated that the function $\Gamma(1-z)I(z,f)$ is analytic on $\mathbb{C}\backslash\{1\}$, while $\Gamma(1-z)$ possesses poles at all integers $n\geq 2$. This establishes that $I(z,f)$ vanishes at these points.
\end{remark}

\subsection{Evaluation of $L(-n,f)$, $n\in \mathbb{Z}_{\geq0}$}\label{subsec5.5}

\begin{theorem}\label{thm5.5}
For all integers $n \geq 0$, the following relation holds: \[L(-n,f)=-\frac{C_{f,n+1}(n_f)}{n+1}=(-1)^n\frac{C_{\underline{f},n+1}(1-n_f)}{n+1}.\]
\end{theorem}

\begin{proof}
Substitute $z$ with $-n$ in Equation (\ref{eqt5.9}), where $n$ is a non-negative integer. We obtain: \[I(-n,f)=\frac{1}{2i\pi}\int_{\mu_{2}}\frac{s^{-n-1}e^{n_{f}s}p_{f}(s)}{1-e^{s}}\diff{s}=\underset{s=0}{\text{Res}}\left(\frac{s^{-n-1}e^{n_{f}s}p_{f}(s)}{1-e^{s}}\right).\] Using (\ref{eqt2.11}) and (\ref{eqt2.2}), we find:
\begin{align*}
\frac{s^{-n-1}e^{n_{f}s}p_{f}(s)}{1-e^{s}}=&-s^{-n-2}e^{n_{f}s}\frac{s}{e^{s}-1}p_{f}(s)\\
=&-s^{-n-2}e^{n_{f}s}f(s)\\
=&-\sum_{p=0}^{+\infty}\frac{C_{f,p}(n_f)}{p!}s^{p-n-2}.
\end{align*}
Therefore, using (\ref{eqt5.11}), we derive:
\[L(-n,f)=\Gamma(1+n)I(-n,f)=-n!\frac{C_{f,n+1}(n_f)}{(n+1)!}=-\frac{C_{f,n+1}(n_f)}{n+1}.\]
Lastly, using (\ref{eqt2.8}), we obtain:
$L(-n,f)=(-1)^n\frac{C_{\underline{f},n+1}(1-n_f)}{n+1}$.
\end{proof}

\subsection{LC-function formula}\label{subsec5.6}

The Hurwitz zeta function satisfies the formula, referred to as Hurwitz's formula, presented in (\cite{Apostol}, Chapter 12, page 256), as follows:
\begin{equation}\label{eqt5.13}
\zeta(1-z,a)=\frac{\Gamma(z)}{(2\pi)^{z}}\left(e^{-\frac{i\pi z}{2}}F(a,z)+e^{\frac{i\pi z}{2}}F(-a,z)\right),
\end{equation}
where $0<a\leq1$, and $F(s,z)$ denotes the periodic zeta function given by: $F(s,z)=\sum_{n=1}^{+\infty}\frac{e^{2in\pi s}}{n^z}$. Prior to generalizing this formula for LC-functions, it is essential to introduce another class of functions, referred to as \textit{FC-functions associated to} $f$. These functions are defined by an integral over a Hankel contour $\mu_{f}$, comprising three parts $\mu_{f,1}$, $\mu_{f,2}$, and $\mu_{f,3}$, as illustrated in Figure 5. 
\begin{figure}[H]
\begin{center}
\includegraphics[]{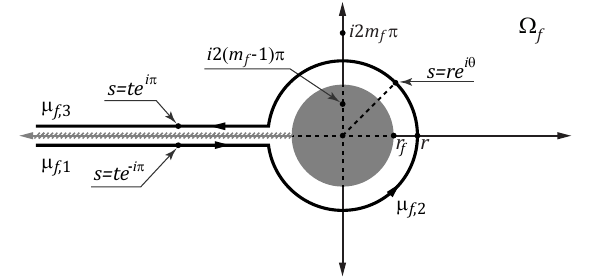}
\caption{The Hankel contour path $\mu_{f}$ is traversed in a counterclockwise direction, maintaining $r_{f}<r<2m_{f}\pi$.
}
\label{fig:Fig5}
\end{center}
\end{figure}
The contour parameters are defined such that $r_{f}<r<2m_{f}\pi$ and $m_{f}$ is the smallest integer satisfying $r_{f}<2m_{f}\pi$. The \textit{FC-functions associated to} $f$ can be represented as:
\begin{equation}\label{eqt5.14}
F(z,f):=\Gamma(1-z)J(z,f),
\end{equation}
where
\begin{equation}\label{eqt5.15}
J(z,f):=\frac{1}{2i\pi}\int_{\mu_{f}}\frac{s^{(z-1,f)}}{e^{-s}-1}\diff{s}.
\end{equation}

\begin{theorem} \label{thm5.6}
The function $J(z,f)$ is entire and the function $F(z,f)$ is analytic for $z\in\mathbb{C}\backslash\{1,2,3...\}$.
\end{theorem}

\begin{proof}
Along $\mu_{f,1}$ and $\mu_{f,3}$, where $s = te^{-i\pi}$ and $s = te^{i\pi}$ respectively, we can infer from Equations (\ref{eqt4.4}) and (\ref{eqt4.11}):\\
$-$ On $\mu_{f,1}$: \[s^{(z,f)}=t^{z}e^{-i\pi z}\sum_{n=0}^{+\infty}\binom{z}{n}P_{f,n}(-t)^{-n}=e^{-i\pi z}\sum_{n=0}^{+\infty}\binom{z}{n}P_{\underline{f},n}t^{z-n}=e^{-i\pi z}t^{(z,\underline{f})}.\]
$-$ On $\mu_{f,3}$: \[s^{(z,f)}=t^{z}e^{i\pi z}\sum_{n=0}^{+\infty}\binom{z}{n}P_{f,n}(-t)^{-n}=e^{i\pi z}t^{(z,\underline{f})}.\]
Therefore, 
\begin{align*}
2i\pi J(z,f)&=\left(\int_{\mu_{f,1}}+\int_{\mu_{f,2}}+\int_{\mu_{f,3}}\right)\frac{s^{(z-1,f)}}{e^{-s}-1}\diff{s}\\
&=e^{-i\pi z}\int_{+\infty}^{r}\frac{t^{(z-1,\underline{f})}}{e^{t}-1}\diff{t}+\int_{\mu_{f,2}}\frac{s^{(z-1,f)}}{e^{-s}-1}\diff{s}+e^{i\pi z}\int_{r}^{+\infty}\frac{t^{(z-1,\underline{f})}}{e^{t}-1}\diff{t}.
\end{align*}
For the integral along $\mu_{f,2}$, a variable substitution $s=re^{i\theta}$ leads to
\begin{align*}
2i\pi J(z,f)&=2i\sin(\pi z)\int_{r}^{+\infty}\frac{t^{(z-1,\underline{f})}}{e^{t}-1}\diff{t}
+ir\int_{-\pi}^{\pi}\frac{e^{i\theta}\left(re^{i\theta}\right)^{(z-1,f)}}{e^{-r e^{i\theta}}-1}\diff{\theta}\\
&=2i\sin(\pi z)\int_{r}^{+\infty}w(z,t)\diff{t}+J_{1}(z,r)\\
&=2i\sin(\pi z)\int_{1}^{+\infty}w(z,t)\diff{t}+2i\sin(\pi z)\int_{r}^{1}w(z,t)\diff{t}+J_{1}(z,r).
\end{align*}

First of all, the functions $s^{(z-1,f)}$ and $s^{(z-1,\underline{f})}$, $(s,z)\in\Omega_{f}\times\mathbb{C}$, are continuous functions of $s$ on $\Omega_{f}$, and as per \Cref{thm4.2}.(i), they are entire functions of $z$. now it remains to prove, by employing holomorphy under the integral sign, that the function $\int_{1}^{+\infty}w(z,t)\diff{t}$ is well-defined and entire. 

From (\ref{eqt4.4}) and (\ref{eqt4.2}) we can write $t^{(z-1,\underline{f})}=S_{p_{\underline{f}}}(t,z-1)$. Within every compact disk $|z|\leq M$ ($M>0$) and for each $t\geq1$, \Cref{lem3.2}.(ii).(a) implies that \[|w(z,t)|=\left\lvert\frac{S_{p_{\underline{f}}}(t,z-1)}{e^{t}-1}\right\rvert\leq\frac{rK_{M}}{r-r_0}\frac{|t^{z-1}|}{e^{t}-1}\leq\frac{rK_{M}}{r-r_0}\frac{t^{M-1}}{e^{t}-1},\] where $r_{0}\in (r_{f},r)$ is an arbitrary real number and the constant $K_M$ is independent of $z$. Given that the integral $\int_{1}^{+\infty}\frac{t^{M-1}}{e^{t}-1}\diff{t}$ converges, holomorphy under the integral sign ensures that $\int_{1}^{+\infty}w(z,t)\diff{t}$ is an entire function.
\end{proof}
We turn now to the LC-function formula stated in the next \Cref{thm5.7}. The proof is written in the same spirit as for the Hurwitz's formula established in (\cite{Apostol}, chapter 12, pages 256-259). First, we present the following lemma which is proved in the same reference.

\begin{lemma}\label{lem5.1}
Let $S(\lambda)$ denote the region that remains when we remove from the s-plane all open circular disks of radius $\lambda$, $0<\lambda<\pi$, with centers at $s=2i n\pi$, $n\in\mathbb{Z}$. Then if $0<a\leq 1$ the function \[\eta(s,a)=\frac{e^{as}}{1-e^s}\] is bounded in $S(\lambda)$. $($the bound depends on $\lambda$.$)$
\end{lemma}

\begin{theorem} [The LC-function formula]\label{thm5.7}
For $\sigma < 0$, the following equation holds:
\begin{equation}\label{eqt5.16}
L\left(1-z,f\right)=\frac{\Gamma(z)}{(2\pi)^{z}}\left(e^{-\frac{i\pi z}{2}}F(z,f_{(-2i\pi)})+e^{\frac{i\pi z}{2}}F(z,f_{(2i\pi)})\right).
\end{equation}
\end{theorem}

\begin{proof}
Let $\omega$ be one of the two values $2i\pi$ or $-2i\pi$ and consider the function \[J_{N}(z,f_{(\omega)})=\frac{1}{2i\pi}\int_{\mu_{f_{(\omega)}}(N)}\frac{s^{\left(z-1,f_{(\omega)}\right)}}{e^{-s}-1}\diff{s},\]
where $\mu_{f_{(\omega)}}(N)$ is the closed loop shown in \Cref{fig:Fig6} with $r_{f_{(\omega)}}<r<2n_{f}\pi$, $N$ being an integer $>n_{f}$. 

\begin{figure}[H]
\begin{center}
\includegraphics[]{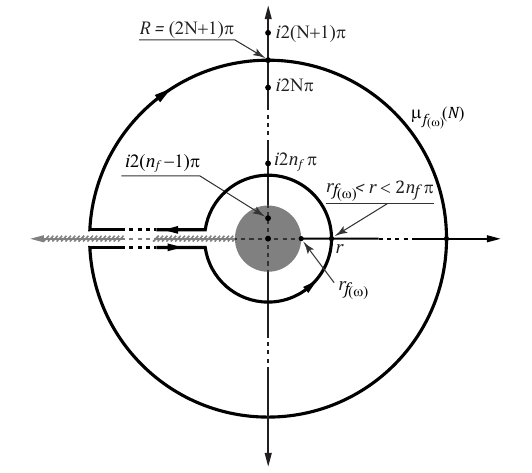}
\caption{Contour plot for deriving the LC-function formula, maintaining $r_{f}<r<2 n_{f}\pi$ and $R=(2N+1)\pi$.}
\label{fig:Fig6}
\end{center}
\end{figure}

From (\ref{eqt4.13}), it follows that $r_{f_{(\omega)}}=2\pi r_{f}$. According to (\ref{eqt5.1}), we have $n_{f}-1\leq r_{f}<n_{f}$, which implies that $2(n_{f}-1)\pi\leq r_{f_{(\omega)}}<2n_{f}\pi$. Thus, $n_f$ represents the smallest integer for which $r_{f_{(\omega)}}<2n_{f}\pi$. As defined in (\ref{eqt5.15}),\[J(z,f_{(\omega)})=\frac{1}{2i\pi}\int_{\mu_{f_{(\omega)}}}\frac{s^{(z-1,f_{(\omega)})}}{e^{-s}-1}\diff{s},\] where $\mu_{f_{(\omega)}}$ denotes the infinite Hankel contour surrounding the negative real axis in the counterclockwise direction, with a radius of $r$ (\Cref{fig:Fig6}).

Suppose $\sigma < 0$. First, we aim to demonstrate that $\lim_{N}J_{N}(z,f_{(\omega)})=J(z,f_{(\omega)})$. In order to establish this, it is sufficient to prove that the integral $J_{R}$ along the outer circle approaches $0$ as $N$ tends to $\infty$. Herein, we express $J_R$ in the following form: \[J_{R}=iR\int_{-\pi}^{\pi}\frac{e^{i\theta}\left(Re^{i\theta}\right)^{(z-1,f_{(\omega)})}}{e^{-R e^{i\theta}}-1}\diff{\theta}=iR\int_{-\pi}^{\pi}e^{i\theta}\eta(R e^{i\theta},1)\left(Re^{i\theta}\right)^{(z-1,f_{(\omega)})}\diff{\theta}.\] Given that \[\left\lvert\left(Re^{i\theta}\right)^{z-1}\right\rvert=R^{\sigma-1}e^{-y\theta}\leq R^{\sigma-1}e^{\pi|y|},\] it follows from (\ref{eqt4.4}), (\ref{eqt4.2}), and \Cref{lem3.2}.(ii).(b) that: \[\left\lvert \left(Re^{i\theta}\right)^{(z-1,f_{(\omega)})}\right\rvert=\left\lvert S_{p_{f_{(\omega)}}}(Re^{i\theta},z-1)\right\rvert\leq\frac{rK_{z}}{r-r_0}\left\lvert \left(Re^{i\theta}\right)^{z-1}\right\rvert\leq\frac{rK_{z}e^{\pi|y|}}{r-r_0}R^{\sigma-1},\] where $r_{0}\in (r_{f_{(\omega)}},r)$ is an arbitrary real number and the constant $K_{z}$ depends only on $z$. According to \Cref{lem5.1}, since the outer circle is contained within the set $S(\lambda)$, $ \left\lvert\eta(R e^{i\theta},1)\right\rvert$ is bounded by a constant $D>0$. Therefore, the integral $J_{R}$ is bounded by: \[\frac{2\pi rDK_{z}e^{\pi|y|}}{r-r_0}R^{\sigma},\] which, given that $\sigma < 0$, approaches $0$ as $R$ tends to $\infty$. Consequently,

\begin{equation}\label{eqt5.17}
\lim_{N}J_{N}(z,f_{(\omega)})=\frac{1}{2i\pi}\int_{\mu_{f_{(\omega)}}}\frac{s^{(z-1,f_{(\omega)})}}{e^{-s}-1}\diff{s}=J(z,f_{(\omega)}).
\end{equation}

We proceed to compute $J_{N}(z,f_{(\omega)})$ explicitly, utilizing Cauchy's residue theorem as follows:

\begin{align*}
J_{N}(z,f_{(\omega)})&=-\sum_{n=n_{f}}^{N}\underset{s=2in\pi}{\text{Res}}\left(\frac{s^{\left(z-1,f_{(\omega)}\right)}}{e^{-s}-1}\right)-\sum_{n=n_{f}}^{N}\underset{s=-2in\pi}{\text{Res}}\left(\frac{s^{\left(z-1,f_{(\omega)}\right)}}{e^{-s}-1}\right)\\
&=\sum_{n=n_{f}}^{N}(2in\pi)^{\left(z-1,f_{(\omega)}\right)}+\sum_{n=n_{f}}^{N}(-2in\pi)^{\left(z-1,f_{(\omega)}\right)}.
\end{align*}

According to (\ref{eqt4.15}) we write \[J_{N}(z,f_{(\omega)})=(2i\pi)^{z-1}\sum_{n=n_{f}}^{N}n^{\left(z-1,f_{(\omega/(2i\pi))}\right)}+(-2i\pi)^{z-1}\sum_{n=n_{f}}^{N}n^{\left(z-1,f_{(-\omega/(2i\pi))}\right)}.\]
Letting $N$ tend to $\infty$ and using (\ref{eqt5.17}) and (\ref{eqt5.2}), given that $\sigma < 0$, we obtain:

\begin{equation}\label{eqt5.18}
J(z,f_{(\omega)})=(2i\pi)^{z-1}L\left(1-z,f_{\left(\frac{\omega}{2i\pi}\right)}\right)+(-2i\pi)^{z-1}L\left(1-z,f_{\left(-\frac{\omega}{2i\pi}\right)}\right).
\end{equation} Thus, by substituting (\ref{eqt5.18}) into (\ref{eqt5.14}), we derive:
\[F(z,f_{(\omega)})=\frac{\Gamma(1-z)}{(2\pi)^{1-z}}\left(i^{z-1}L\left(1-z,f_{\left(\frac{\omega}{2i\pi}\right)}\right)+i^{1-z}L\left(1-z,f_{\left(-\frac{\omega}{2i\pi}\right)}\right)\right).\]
Upon replacing $\omega$ with $2i\pi$ and subsequently with $-2i\pi$, and noting that $f_{(-1)}=\underline{f}$ and $f_{(1)}=f$, we obtain:
\[F(z,f_{(2i\pi)})=\frac{\Gamma(1-z)}{(2\pi)^{1-z}}\left(i^{z-1}L\left(1-z,f\right)+i^{1-z}L\left(1-z,\underline{f}\right)\right)\] and \[F(z,f_{(-2i\pi)})=\frac{\Gamma(1-z)}{(2\pi)^{1-z}}\left(i^{z-1}L\left(1-z,\underline{f}\right)+i^{1-z}L\left(1-z,f\right)\right).\]
Besides, $e^{i\pi z/2}=i^{z}$ and $e^{-i\pi z/2}=i^{-z}$. Hence \[e^{\frac{i\pi z}{2}}F(z,f_{(2i\pi)})=\frac{\Gamma(1-z)}{i(2\pi)^{1-z}}\left(e^{i\pi z}L\left(1-z,f\right)-L\left(1-z,\underline{f}\right)\right)\] and
\[e^{\frac{-i\pi z}{2}}F(z,f_{(-2i\pi)})=\frac{\Gamma(1-z)}{i(2\pi)^{1-z}}\left(L\left(1-z,\underline{f}\right)-e^{-i\pi z}L\left(1-z,f\right)\right).\] Thus, the formula (\ref{eqt5.16}) is ultimately derived by summing the two equations above and employing the complementary formula $\Gamma(z)\Gamma(1-z)=\frac{\pi}{\sin(\pi z)}=\frac{2i\pi}{e^{i\pi z}-e^{-i\pi z}}$.

\end{proof}

\newpage

\bibliographystyle{amsplain}
\providecommand{\bysame}{\leavevmode\hbox to3em{\hrulefill}\thinspace}

\end{document}